\documentclass[entropy,article,accept,pdftex,oneauthor]{Definitions/mdpi} 
\usepackage[utf8]{inputenc}
\usepackage{amsmath}
\usepackage{calc}
\usepackage{accents}
\usepackage[symbol]{footmisc}
\newcommand{\dbtilde}[1]{\accentset{\approx}{#1}}

\firstpage{1} 
\pubvolume{1}
\issuenum{1}
\articlenumber{0}
\pubyear{2025}
\copyrightyear{2025}
\datereceived{ } 
\daterevised{ } 
\dateaccepted{ } 
\datepublished{ } 
\hreflink{https://doi.org/} 

\Title{NUMERICAL INTEGRATION OF STOCHASTIC DIFFERENTIAL EQUATIONS: THE HEUN ALGORITHM REVISITED AND THE IT\^O STRATONOVICH CALCULUS}

\TitleCitation{Numerical Integration of SDEs: The Heun Algorithm Revisited and Itô-Stratonovich Calculus}

\Author{Riccardo Mannella $^{1}$\orcidA{}}

\AuthorNames{Riccardo Mannella}

\isAPAStyle{%
        \AuthorCitation{Riccardo Mannella}
          }{%
         \isChicagoStyle{%
         \AuthorCitation{Riccardo Mannella}
         }{
         \AuthorCitation{Riccardo Mannella}
         }
}

\address{%
$^{1}$ \quad Dipartimento di Fisica, Università di Pisa, Italy

\corres{Correspondence: riccardo.mannella@unipi.it}}

\abstract{
The widely used Heun algorithm for the numerical integration of stochastic differential equations (SDEs) is critically re-examined. We discuss and evaluate several alternative implementations, motivated by the fact that the standard Heun scheme is constructed from a low-order integrator. The convergence, stability, and equilibrium properties of these alternatives are assessed through extensive numerical simulations. 
Our results confirm that the standard Heun scheme remains a benchmark integration algorithm for SDEs due to its robust performance. As a byproduct of this analysis, we also disprove a previous claim in the literature regarding the strong convergence of the Heun scheme.
}


\begin{document}


\section{Introduction}
Stochastic Differential Equations (SDEs) are an indispensable tool for modeling complex systems across a vast range of scientific disciplines, including physics, chemistry, biology, and finance~\citep{Gardiner_2010,Oksendal_2003}. These equations provide a powerful framework for describing systems that are subject to random fluctuations, or "noise," which is often an intrinsic feature of the underlying dynamics. Unlike Ordinary Differential Equations (ODEs), which have deterministic solutions, the solution to an SDE is itself a stochastic process, representing a collection of possible trajectories.

A central challenge in working with SDEs is that most non-trivial cases cannot be solved analytically. Consequently, their analysis relies heavily on numerical integration schemes to approximate solutions~\citep{Platen_Kloeden_1992}. The development of robust and efficient numerical algorithms is therefore of paramount importance. An ideal scheme must balance several competing factors: it should be accurate, converging to the true solution as the integration time step decreases; it must be stable, avoiding divergence even for large time steps or stiff problems; and it should be computationally efficient.

Among the many algorithms developed for this purpose, the Heun scheme has emerged as a widely used and reliable method. As a stochastic extension of the classic Runge-Kutta predictor-corrector method, it is known for its favorable balance of simplicity, stability, and accuracy. Its structure naturally handles the Itō-Stratonovich dilemma, converging to the Stratonovich interpretation of the SDE, which is often preferred in physical modeling where the noise represents a smoothed, physical process rather than a purely mathematical construct~\citep{vanKampen_2007,Mannella_McClintock_2012}.

Despite its widespread use, the standard Heun scheme is built upon the Euler-Maruyama method, a scheme with a relatively low order of strong convergence. This fact motivates a critical re-examination: Could the performance of the Heun method be improved by incorporating higher-order building blocks? How do its stability and accuracy compare to other, more complex Taylor-based or Runge-Kutta-type schemes?

This paper addresses these questions through a systematic investigation of the Heun algorithm and its variants. We critically evaluate several alternative implementations, comparing their performance in terms of strong convergence, numerical stability, and their ability to reproduce the correct long-time equilibrium distributions of a non-trivial SDE. Our findings confirm that the standard Heun scheme remains a benchmark algorithm for its robustness and reliability. As a byproduct of this analysis, we also re-evaluate and ultimately disprove a previous claim in the literature regarding the strong convergence order of the Heun scheme.


\section{The Heun Algorithm and Itô-Stratonovich Calculus}

The \textit{Heun algorithm} for the numerical integration of Stochastic Differential Equations (SDEs) is a stochastic extension of an integrator first introduced by Heun in 1900~\citep{HEUN1900}. That paper addressed the then-nascent field of Runge-Kutta schemes for integrating Ordinary Differential Equations (ODEs)~\citep{BUTCHER1996247}.

Given the ODE:
\begin{equation}
    \dot{x} = f(x,t),
    \label{ODE1}
\end{equation}
the Heun scheme calculates the next step, $x_{n+1} \equiv x(t+h)$, from the current step, $x_n \equiv x(t)$, as follows:
\begin{align}
    \tilde{x}   &= x_n + h f(x_n,t) \nonumber \\
    \dbtilde{x} &= x_n + h f(\tilde{x},t+h) \nonumber \\
    x_{n+1}     &= \frac{1}{2} (\tilde{x} + \dbtilde{x}) = x_n + \frac{h}{2} [f(x_n,t) + f(\tilde{x},t+h)].
    \label{Heun1}
\end{align}
Equation~\ref{Heun1} shows that the Heun scheme is a \textit{predictor-corrector} method. It first \textit{predicts} a value with an explicit Euler step and then \textit{corrects} it by averaging this explicit step with an implicit Euler step.

One of the earliest papers to mention a stochastic version of the Heun algorithm is~\citep{Rumelin_1982}, where this name was used to refer to a stochastic Runge-Kutta scheme derived in~\citep{Kloeden_Pearson_1977}. Given the one-dimensional SDE:
\begin{equation}
    dx = f(x,t) dt + g(x,t) dW,
    \label{SDE1}
\end{equation}
where $W$ is a standard Wiener process, the well-known Euler-Maruyama (in the following, Euler) scheme in the \textit{It\^o sense} is:
\begin{equation}
    x_{n+1} = x_n + h f(x_n,t) + g(x_n,t) Z_1,
    \label{SDEEuler}
\end{equation}
where $Z_1 = Z_1(h) \equiv \int_t^{t+h} dW = \sqrt{h} Y_1$, with $Y_1 \sim \mathcal{N}(0, 1)$, is a Gaussian random variable with zero mean and variance $h$. In~\citep{Rumelin_1982}, the Heun scheme for SDEs was written as:
\begin{equation}
    x_{n+1} = x_n + \frac{h}{2} [f(x_n,t) + f(\tilde{x},t+h)] + \frac{1}{2} [g(x_n,t) + g(\tilde{x},t+h)] Z_1,
    \label{SDEHeun}
\end{equation}
where the predictor step is:
\begin{equation}
    \tilde{x} = x_n + h f(x_n,t) + g(x_n,t) Z_1, 
    \label{predictor}
\end{equation}
and the same random variable $Z_1$ is used in both stages. Much like its deterministic counterpart, the stochastic Heun scheme can be seen as the average of an explicit and an implicit stochastic Euler step.

It is well known (see, for example,~\citep{Rumelin_1982,Kloeden_Pearson_1977,Klauder_Petersen_1985}) that the stochastic evolution described by Eq.~\ref{SDEHeun} converges to the evolution of Eq.~\ref{SDE1} interpreted in the \textit{Stratonovich} sense. This is equivalent to the following SDE in the Itô sense:
\begin{equation}
    dx = \left(f(x,t) + \frac{1}{2} g(x,t) \frac{d g(x,t)}{dx}\right) dt + g(x,t) dW.
    \label{SDE2}
\end{equation}

It was claimed in~\citep{Kloeden_Pearson_1977} that the Heun scheme is a scheme with a strong convergence order $O(h^{3/2})$\footnote{See Section~\ref{section:convergence} for a definition of strong convergence. Also, we will loosely use $O(h^\alpha)$ and $o(h^\alpha)$ when
terms proportional to $h^\alpha$ are kept or discarded, respectively.}
Furthermore, in~\citep{Mannella_2002}, it was demonstrated that for a potential $V(x)$ bounded from below (i.e., $\inf_{x} V(x) > -\infty$), with $f(x) = -V'(x)$ and $g(x) = \sqrt{2D}$, the scheme in Eq.~\ref{SDEHeun} yields the correct equilibrium distribution for $x$ up to order $o(h^2)$.

In the following sections, we will consider only single-step algorithms with a time step $h$. For simplicity, we will use the notation $x(h) \equiv x_{n+1}$ and $x(0) \equiv x_n$. When there is no ambiguity, we will use shorthand such as $f \equiv f(x,t)$ and $g' \equiv \frac{d g(x,t)}{dx}$. Subscripts will be added when necessary to avoid confusion (e.g., $f_0 \equiv f(x_0,t_0)$).

Although it may seem unsurprising that the Heun scheme yields a Stratonovich evolution—given that its elementary step averages quantities evaluated at the beginning and end of the time interval, which mirrors the definition of Stratonovich calculus~\citep{Mannella_McClintock_2012}—several questions naturally arise:
\begin{itemize}
    \item The Heun scheme produces a Stratonovich evolution by combining two Euler schemes that use the Itô prescription. Should one instead use Euler schemes derived from the Stratonovich prescription?
    \item The standard stochastic Euler scheme has a strong convergence order of $O(h^{1/2})$~\citep{Platen_Kloeden_1992}, unlike the deterministic Euler scheme. Would it be beneficial to replace it in the Heun scheme with an integrator that has a strong convergence order of $O(h)$?
    \item How does the Heun scheme compare with other Taylor-based higher-order schemes?
\end{itemize}
To address these points, some background is necessary (see also~\citep{Mannella_McClintock_2012}).

\subsection{The Stochastic Euler Scheme and Stochastic Calculus}
An elementary $O(h)$ scheme for Eq.~\ref{SDE1}, obtained via a Taylor expansion around $x_0$, is given by (see, for example,~\citep{Mannella_1989}):
\begin{equation}
    x(h) = x(0) + g_0 Z_1 + h f_0 + g_0 g'_0 \int_0^h W(s) \,dW(s).
    \label{SDEEuleroWithIntegral}
\end{equation}
If Itô calculus is used, the integral is $\int_0^h W(s) \,dW(s) = \frac{1}{2} (Z_1^2 - h)$, which leads to the Milstein scheme (see Eq.~\ref{milstein} below). If this term is approximated as zero, the standard Euler scheme (Eq.~\ref{SDEEuler}) is recovered. If Stratonovich calculus is used, the integral is $\int_0^h W(s) \circ dW(s) = \frac{1}{2} Z_1^2$, which yields the Euler-Stratonovich scheme:
\begin{equation}
    x(h) = x(0) + g_0 Z_1 + h f_0 + \frac{1}{2} g_0 g'_0 Z_1^2.
    \label{SDEEulerStratonovich}
\end{equation}
The Heun algorithm requires an Euler step involving the final point, $x(h)$. For It\^o calculus, this step corresponds to Eq.~\ref{SDEEuler}, replacing $x_n$ by $x_{n+1}$ in the arguments of $f(x)$ and $g(x)$. In contrast, when using Stratonovich calculus, expanding around $x(h)$ results in:
\begin{equation}
    x(h) = x(0) + g_h Z_1 + h f_h - \frac{1}{2} g_h g'_h Z_1^2.
    \label{SDEEulerStratonovichFinal}
\end{equation}

\subsection{Strong Approximation $O(h)$: The Milstein Scheme}
A strong approximation $O(h)$  for integrating Eq.~\ref{SDE1} was proposed by Milstein~\citep{Milshtejn_1975}:
\begin{equation}
    x(h) = x(0) + g_0 Z_1 + h f_0 + \frac{1}{2} g_0 g'_0 (Z_1^2 - h).
    \label{milstein}
\end{equation}
Comparing this to Eq.~\ref{SDE2}, the Milstein scheme corresponds to applying the Euler-Stratonovich scheme (Eq.~\ref{SDEEulerStratonovich}) to the following SDE\footnote[2]{The symbol $\circ$ indicates that this SDE is to be interpreted in the Stratonovich sense.}:
\begin{equation}
    dx = \left(f(x,t) - \frac{1}{2} g(x,t) \frac{d g(x,t)}{dx}\right) dt + g(x,t) \circ dW.
    \label{SDE3}
\end{equation}

\subsection{Higher-Order Taylor-based Schemes}
A Taylor-based scheme of order $O(h^2)$ using Itô calculus was derived in~\citep{Rao_1974}, and a revised version using Stratonovich calculus was developed in~\citep{Mannella_1989}. Even higher-order Taylor-based schemes were derived in~\citep{Platen_Kloeden_1992}, leveraging a chain rule for higher-order multiple stochastic integrals.

As shown both analytically and numerically in~\citep{Mannella_2002}, care must be taken when deriving higher-order schemes. Specifically, if the terms containing stochastic integrals are computed to a certain order, the deterministic terms must be kept at the same order. One should not include higher-order deterministic terms, even if they are easy to derive\footnote{The deterministic terms are often simpler to compute, which might tempt one to include them regardless of the order achieved in the stochastic part.}. Doing so can result in a scheme that performs worse than one where all terms—both stochastic and deterministic—are consistently maintained at the same order of accuracy.

In the following, we will use the results of~\citep{Mannella_1989,Mannella1989FastAP}, which are reproduced here for convenience:
\begin{flalign}
x(h) = & x(0) + g_0 Z_1 + h f_0 + \frac{1}{2} g_0 g_0' Z_1^2 && O(h^{1/2})+O(h) \nonumber \\
       & +\left( g_0 f_0' - f_0 g_0' \right) Z_2 + h f_0 g_0' Z_1
       + \frac{1}{3!} g_0 \left( g_0' g_0'+ g_0'' g_0 \right) Z_1^3 && O(h^{3/2}) \nonumber \\
       & + \frac{h^2}{2} f_0 f_0' + \frac{1}{2} f_0' g_0' g_0 \left( Z_1 Z_2 - Z_3\right) 
       + \frac{1}{2} f_0'' g_0 g_0 \left( Z_1 Z_2 - Z_3\right) && O(h^2)\nonumber \\
       & + g_0' \left( g_0 f_0' - f_0 g_0' \right) Z_3 + 
       g_0'g_0' f_0 \frac{1}{2} \left( h Z_1^2 - \frac{1}{2}Z_1 Z_2+ \frac{1}{2}Z_3 \right) && O(h^2)\nonumber \\
       &+ \frac{1}{4!} g_0 \left(g_0^2 g_0''' + g_0' (g_0 g_0')'\right) Z_1^4 +\frac{h}{4} g_0'' g_0 f_0 \left( Z_1^2 - \frac{1}{2}(Z_1 Z_2 - Z_3)\right) && O(h^2)
       \label{CUP}
\end{flalign}
where two additional stochastic integrals appear: $Z_2 = Z_2(h) \equiv \int_0^h Z_1(s) ds = \int_0^h \left(\int_0^s dW \right) ds$ and $Z_3 = Z_3(h) \equiv \int_0^h Z_2(s) dW(s)$. A suitable representation for $Z_2$ is easily found, as it is Gaussian~\citep{Rao_1974, Mannella_1989}:
\begin{equation}
    Z_2 = Z_2(h) = h \left(\frac{Z_1}{2} + \frac{\sqrt{h}}{2 \sqrt{3}} Y_2\right),
\end{equation}
where $Y_2 \sim \mathcal{N}(0, 1)$ is independent of $Y_1$. Finding a representation for $Z_3$ is much harder, as a number of constraints must be satisfied. Considering the lowest orders in $h$, the following moments are found:
\begin{gather*}
    \langle Z_3(h) \rangle = 0, \quad \langle [Z_3(h)]^2 \rangle = \frac{h^4}{12} \quad
    \langle Z_1(h)Z_3(h) \rangle = 0, \quad \langle Z_2(h)Z_3(h) \rangle = 0 \\
    \langle [Z_1(h)]^2 Z_3(h) \rangle = \frac{1}{4}h^3, \quad \langle Z_1(h)Z_2(h)Z_3(h) \rangle = \frac{h^4}{12} \quad
    \langle [Z_2(h)]^2 Z_3(h) \rangle = \frac{h^5}{15},
\end{gather*}
where $\langle \dots \rangle$ denotes an average over stochastic realizations. The number of conditions to satisfy exceeds the number of unknown quantities, so an arbitrary choice must be made. In~\citep{Mannella_1989}, the representation
\begin{equation}
    Z_3(h) = \frac{h^2}{6} \{ Y_1^2 - h +Y_3 \}
    \label{z3cup}
\end{equation} 
was used\footnote{A misprint in Ref~\citep{Mannella_1989} has been corrected here.}, where $Y_3 \sim \mathcal{N}(0, 1)$ is independent of $Y_1$ and $Y_2$. In preparation for this manuscript, an
alternative representation was derived:
\begin{equation}
    Z_3 =  Z_1 Z_2 - \frac{h^2}{2} \left( 1 + \frac{Y_3}{3}\right).
    \label{z3gemini}
\end{equation}

\section{Algorithms and Dynamical System}

The numerical schemes studied here are typically used to integrate SDEs to derive equilibrium or quasi-equilibrium properties. These include calculating the mean first passage time to a boundary, deriving reaction rate constants, and determining the long-time equilibrium distribution.

It is therefore relevant to test these approaches using a non-trivial SDE with a known equilibrium distribution. Verifying agreement with the full distribution is a more stringent test than merely checking the behavior of a few moments or cumulants. This focus on long-time accuracy is at odds with the typical analysis of integration schemes, which often investigates the behavior of stochastic trajectories over short time scales.

However, it is possible in some cases to connect these two regimes. Using the formalism of~\citep{Sancho_Miguel_Katz_Gunton1982}, one can infer the equilibrium distribution generated by a given numerical scheme from its short-time properties~\citep{Mannella_2002, Mannella_2004}. This connection demonstrates that schemes with more accurate short-time behavior also produce numerical equilibrium distributions that are closer to the theoretical one. 

The dynamical system used for the test is one of the simplest non-trivial models:
\begin{equation}
    dx = -x (1 + x^2) dt + \sqrt{2 D} (1 + x^2 ) dW = f(x) dt + g(x) dW.
    \label{model}
\end{equation}
This model has the equilibrium distribution
\begin{equation}
P_{eq}(x) = \frac{N}{(1 + x^2)^{1+\alpha+n}},
\label{equilibrium}
\end{equation}
where $\alpha = 1/(2 D)$, $n$ is 0 for Stratonovich calculus and 1 for Itô calculus, and
$N = \Gamma(1+n+\alpha)/(\Gamma(3/2+\alpha)\sqrt{\pi})$ is a normalization constant. It is noteworthy that the presence of a power-law tail makes the comparison with the numerical schemes particularly interesting, because the different schemes must be able to integrate correctly even when the variable $x$ becomes large.

The algorithms tested are as follows, where the labels will be used for identification below:
\begin{description}
    \item[Euler:] The standard Euler scheme, Eq.~\ref{SDEEuler}.
    \item[Heun:] The standard Heun scheme, Eq.~\ref{SDEHeun}.
    \item[Stra:] The Euler-Stratonovich scheme, Eq.~\ref{SDEEulerStratonovich}.
    \item[Miln:] A modified Heun scheme where the Milstein algorithm (Eq.~\ref{milstein}) is used as the basic block for both the predictor and corrector steps: in the corrector step, $f(x)$ and $g(x)$ are evaluated at the $x$ found in the predictor step.
    \item[HeSt:] A modified Heun scheme using the Euler-Stratonovich scheme of Eq.~\ref{SDEEulerStratonovich} for the predictor step and Eq.~\ref{SDEEulerStratonovichFinal} for the corrector step.
    \item[HePC:] An iterated Heun scheme. First, Eq.~\ref{SDEHeun} is used to obtain a tentative final point. This point is then re-inserted into Eq.~\ref{SDEHeun} in place of $\tilde x$. In our simulations, this loop was iterated four times.
    \item[Mil-:] A modified Heun scheme similar to \textbf{Miln} where for the corrector step the $-$ sign in front of the term $\frac{1}{2} g g' $ is used (this resembles what is done in Eq.~\ref{SDEEulerStratonovichFinal} for the Euler-Stratonovich
    scheme)
    \item[T3/2:] The scheme of Eq.~\ref{CUP} including terms up to the ones marked $O(h^{3/2})$
    \item[HPC-:] An iterated Heun scheme similar to \textbf{HePC}, but
    using Stratonovich schemes Eqs~\ref{SDEEulerStratonovich} 
    and~\ref{SDEEulerStratonovichFinal} rather than Euler. 
    \item[CUP1:] The scheme of Eq.~\ref{CUP} using Eq.~\ref{z3cup}
    \item[CUP2:] The scheme of Eq.~\ref{CUP} using Eq.~\ref{z3gemini}
    \item[RK:] The efficient RK scheme of~\citep{Bogoi_et_al_2023}
\end{description}

\section{Numerical results and discussion}

All codes were written in Fortran 90 and parallelized using the MPICH~\citep{mpich} library on machines with 24, 30, or 40 cores. The random number generator used is RAN2~\citep{press_etal:1992} from Numerical Recipes. 

\subsection{Convergence\label{section:convergence}}

We assessed the convergence rate of the different schemes following the approach of~\citep{Higham_2001}, where further details can be found. The core idea is that if $x_f$ is the final point of a trajectory integrated with a time step $h$ from an initial point $x_0$, and $x_{limit}$ is the "true" final point obtained in the limit of $h \rightarrow 0$ starting from the same $x_0$, the expected error satisfies $\mathbb{E} |x_f-x_{limit}| \le A h^\alpha$ for an algorithm with a strong convergence of order $\alpha$.

In the case studied in~\citep{Higham_2001}, $x_{limit}$ was found by solving the SDE analytically. However, since the exact solution to the SDE in Eq.~\ref{model} is not known, we adopted the following procedure for each algorithm:
\begin{itemize}
\item A random initial point was chosen from the Stratonovich equilibrium 
distribution~\ref{equilibrium}.
\item The trajectory was integrated from the initial point up to time $t=1.0$ with a very small time step $h_s=1/2^{17} \approx 7.63 \times 10^{-6}$, storing the random variates generated at each step. This trajectory serves as the "reference" trajectory, and its final point $x_{ref}$ was stored.
\item The integration was repeated from the same initial point, but with a larger integration 
time step $h_n = 2^n h_s$, for $1 \le n \le 13$, and the final point $x_f(h_n)$ was stored. To generate the noise for these larger time steps, the random increments from the reference trajectory were summed. For example, if the noise increments for the reference trajectory from $t$ to $t+h_s$ and from $t+h_s$ to $t+2 h_s$ were $w_1$ and $w_2$, respectively, the noise increment for a trajectory with time step $h_1 = 2 h_s$ from $t$ to $t+2 h_s$ was taken as $w_1+w_2$. This procedure was extended straightforwardly for cases requiring multiple random terms.
\item The absolute error $| x_f(h_n) - x_{ref} |$ was computed and stored.
\item The procedure was repeated for $N=4 \times 10^6$ trajectories, and the average error, $\mathbb{E}[| x_f(h_n) - x_{ref} |]$, was computed.
\item Finally, a fit of the average error versus $h$ was performed using MINUIT2~\citep{James:1975dr} to determine the parameters $A$ and $\alpha$.
\end{itemize}
For integration time steps where $h \gg h_{s}$, it is reasonable to consider $x_{ref}$ a proxy for the true value $x_{limit}$. The fit was performed using simulation data in the range $10^{-4} \le h \le 10^{-2}$.

\begin{figure}[H]
\includegraphics[width=0.8 \textwidth]{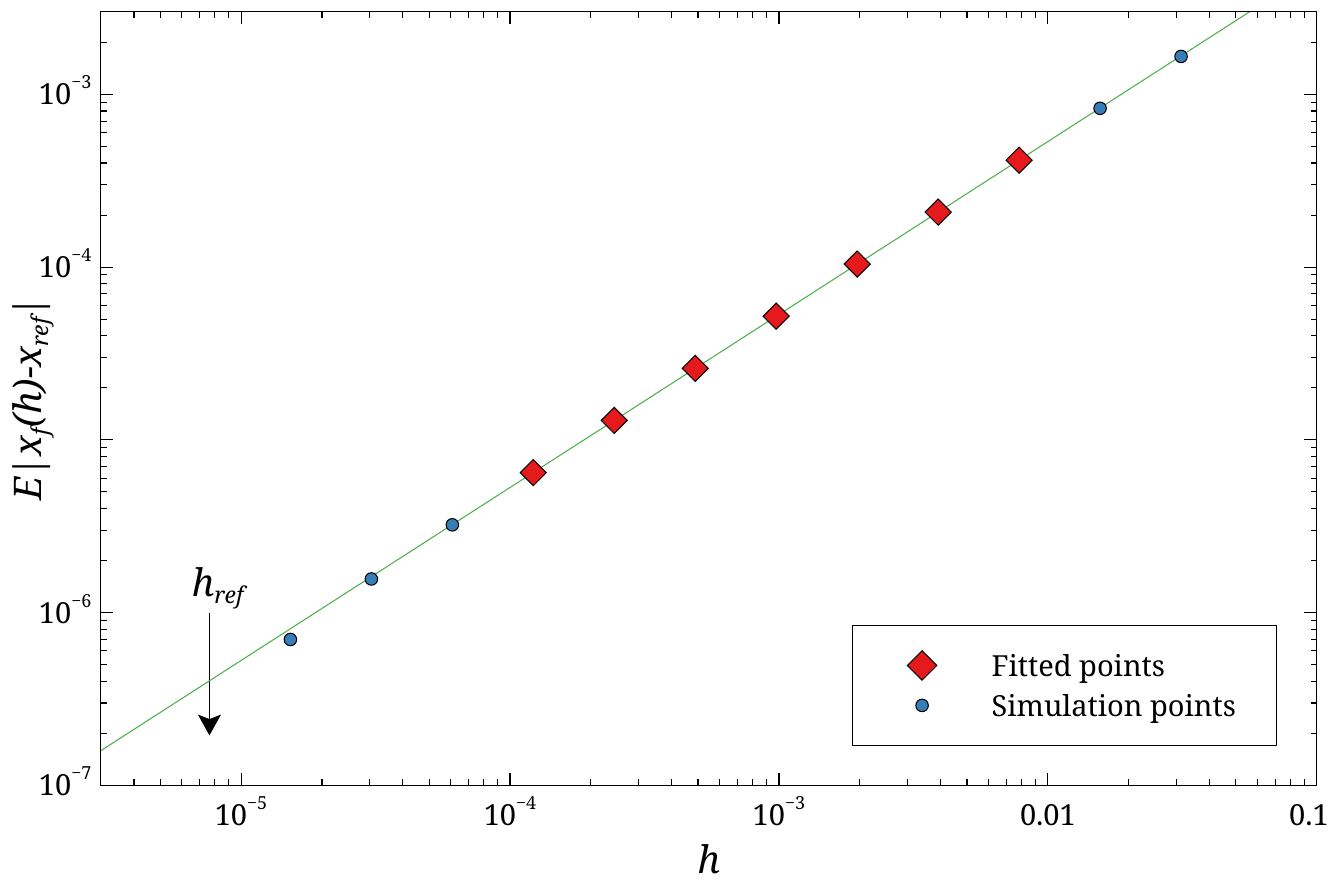}
\caption{Typical data used for the convergence fit, showing
$\mathbb{E} |x_f(h)-x_{ref}|$ vs $h$. The data shown are for the Heun algorithm. The solid line is the fit to the model $f(h)=A h^\alpha$. The symbols are the results of simulations: the red diamonds were used for the fit, while the light blue circles were excluded. The arrow labeled $h_{ref}$ marks the $h$ used for the reference trajectory. The statistical error due to the finite number of averages is smaller than the symbol size. \label{figureconvergency}}
\end{figure}    

A typical case is shown in Fig.~\ref{figureconvergency}. Note that the values of $h$ used in the fit are significantly larger than $h_{ref}$. The results of the fit are summarized in Table~\ref{convergency}. The points marked with light blue circles were excluded from the fit but are shown to verify the expected linear relationship on a log-log plot within the fitting range. To validate our method, we first tested the standard Euler scheme. As expected, the fit yields $\alpha = 0.53 \pm 0.02$, which is consistent with the theoretical convergence order of 0.5 for the Euler scheme.

\begin{table}[H] 
\caption{Table summarizing strong convergence and stability for each scheme.
The convergence order is in the form $\propto A h^\alpha$, and the table shows the value
of the parameters $A$ and $\alpha$ found in the fit; in the simulations for the convergence, $D=0.05$. The column "Stability" reports the number of red dots appearing on each plate of Fig.~\protect{\ref{figurestability}}.\label{convergency}}
\begin{tabularx}{\textwidth}{CCCC}
\toprule
\textbf{Algorithm} & $A$ & $\alpha$ & Stability\\
\midrule
Euler		& $(1.8 \pm 0.3) \times 10^{-2}$	 &  $0.53 \pm 0.02$ &\\
Heun		& $(5.4 \pm 0.5) \times 10^{-2}$     &  $1.00 \pm 0.01$ & 191\\
Stra		& $(9.6 \pm 0.9) \times 10^{-2}$	 &  $1.01 \pm 0.01$ & 115\\
Miln		& $(1.8 \pm 0.2) \times 10^{-2}$	 &  $0.52 \pm 0.01$ &  91\\
HeSt		& $(5.9 \pm 0.6) \times 10^{-2}$	 &  $1.00 \pm 0.01$ & 172\\
HePC		& $(5.4 \pm 0.5) \times 10^{-2}$	 &  $1.00 \pm 0.01$ & 169\\
Mil-		& $(5.6 \pm 0.5) \times 10^{-2}$	 &  $1.00 \pm 0.01$ & 165\\
HPC-		& $(6.2 \pm 0.6) \times 10^{-2}$	 &  $1.00 \pm 0.01$ & 132\\
T3/2    	& $(9.2 \pm 0.9) \times 10^{-2}$	 &  $1.01 \pm 0.01$ & 134\\
CUP1       	& $(7.7 \pm 0.7) \times 10^{-2}$	 &  $1.01 \pm 0.01$ & 117\\
CUP2       	& $(7.7 \pm 0.7) \times 10^{-2}$	 &  $1.01 \pm 0.01$ & 117\\
RK         	& $(2.4 \pm 0.2) \times 10^{-2}$	 &  $0.85 \pm 0.01$ & 135\\
\bottomrule
\end{tabularx}
\end{table}
Most schemes show strong convergence of order $O(h)$. The simple Stratonovich scheme (\textbf{Stra}, Eq.~\ref{SDEEulerStratonovich}) also has $O(h)$ convergence. This may seem surprising given that the Euler scheme is only $O(h^{1/2})$, but it is consistent with the fact that the Stratonovich scheme includes terms up to $O(h)$ from the Taylor expansion in Eq.~\ref{CUP}.

Interestingly, the Heun scheme variant that uses the Milstein algorithm (\textbf{Miln}) as its building block achieves only $O(h^{1/2})$ strong convergence. However, $O(h)$ convergence is recovered in the \textbf{Mil-} scheme, where the sign of the $\frac{1}{2} g g'$ term is changed in the corrector step, a modification made in the spirit of the derivation of Eq.~\ref{SDEEulerStratonovichFinal}. Note, however, that the Heun scheme as a strong convergence $O(h)$, contrary to the claim of~\citep{Kloeden_Pearson_1977}: a similar result was found in~\citep{Bogoi_et_al_2023}

On the other hand, the strong convergence found for the higher-order Taylor-based  schemes is rather surprising.
The \textbf{T3/2} algorithm shows strong convergence of only $O(h)$, despite being based on a Taylor expansion that includes terms up to $O(h^{3/2})$. This issue is even more acute for the \textbf{CUP1} and \textbf{CUP2} schemes, which also exhibit strong convergence of only $O(h)$, despite their Taylor expansion (Eq.~\ref{CUP}) including terms up to $O(h^2)$. We also examined
the Runge-Kutta scheme of~\citep{Honeycutt_1992},
which includes terms $O(h^2)$ and was believed to have a strong convergence $O(h^{3/2})$, and found its strong convergence to be only $O(h)$. To rule out possible numerical artifacts or coding errors, we repeated the assessment of the
Taylor-based schemes using SDE $dx = x dt + D \, x \circ dW$ which can be solved analytically: it was possible to confirm that the integration routines were integrating correctly and the overall numerical approach was sound.

\begin{figure}[H]
\isPreprints{\centering}{} 
\includegraphics[width=\textwidth]{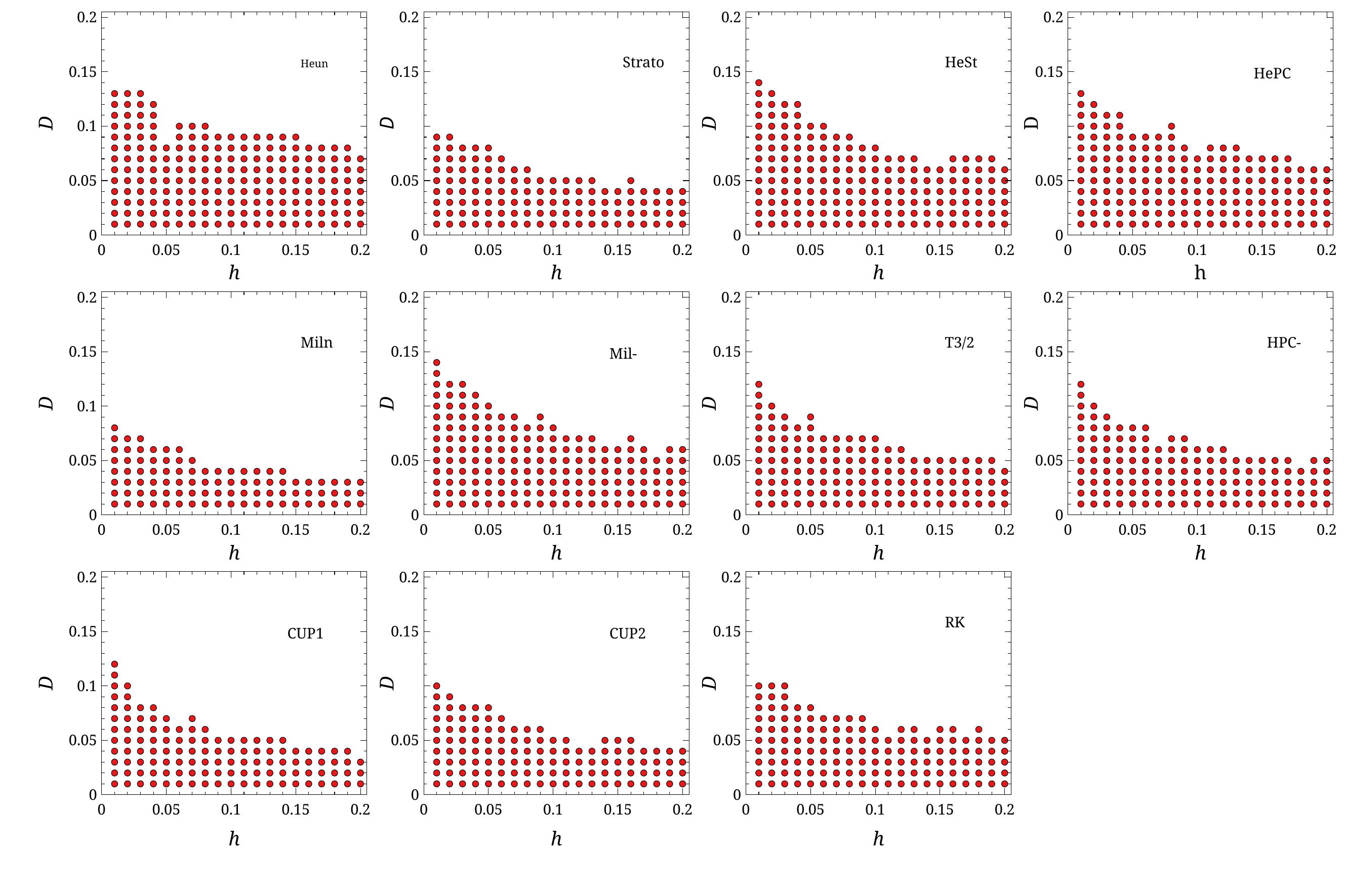}
\caption{Stability regions for the different algorithms, as function
of $D$ and $h$. The red dots mark the parameters for which the given scheme did not blow up.\label{figurestability}}
\end{figure} 

\subsection{Stability}
To establish the parameter regions where the integration remains stable, the system was simulated by varying the parameters $D$ and $h$: 20 values for $D$ and 20 values for $h$ were selected, as shown in Fig.~\ref{figurestability}.
For all simulations and algorithms, $3 \times 10^4$ trajectories, each starting from a random point drawn from the equilibrium distribution, were tracked up to a time of $10^4$. Thus, if $h=10^{-2}$, each trajectory was integrated for $10^4/h = 10^6$ elementary time steps. If $|x|>100$ at any point during the integration for any of the trajectories, we marked that parameter set as unstable for the given algorithm.

The results of these simulations are shown in Fig.~\ref{figurestability}. The "Stability" column in Table~\ref{convergency} reports the number of $(D, h)$ parameter pairs for which the simulations remained stable (i.e., did not diverge) for each algorithm.

The \textbf{Heun} scheme proved to be the most stable, followed by \textbf{HeSt}, \textbf{HePC}, and \textbf{Mil-}. Conversely, the Taylor-based schemes performed worse, particularly \textbf{CUP1} and \textbf{CUP2}. This result is not surprising, as the introduction of a corrector step in a numerical integration scheme is known to improve stability.

\subsection{Stationary distributions}

To study the stationary distributions, we performed simulations across a grid of parameters, using 20 values for $h$ in the range $[0.01,0.2]$ and 10 values for $D$ in the range $[0.01,0.25]$. 
For each integration scheme and for each $(D, h)$ pair, we computed $4 \times 10^5$ trajectories with initial conditions drawn from the Stratonovich theoretical equilibrium distribution. Each trajectory was integrated up to a time $T=10^4$ and sampled at intervals of $\Delta T=1.0$. The value of $x$ at each sample was accumulated to construct the equilibrium probability density, hereafter denoted as $P(x)$. In total, each $P(x)$ is constructed from $4 \times 10^9$ data points.


\begin{figure}[H]
\isPreprints{\centering}{} 
\includegraphics[width=0.48 \textwidth]{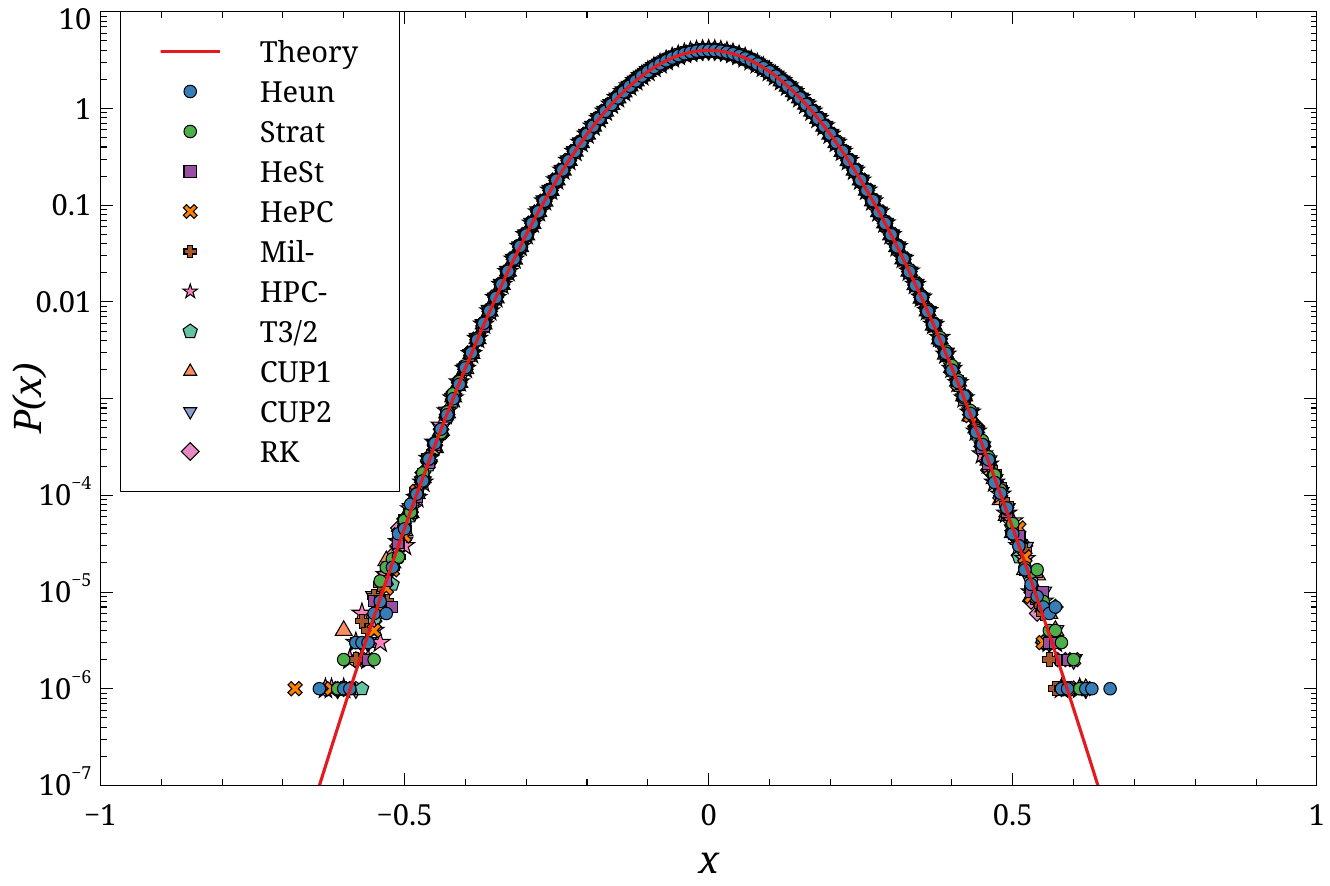}
\includegraphics[width=0.48 \textwidth]{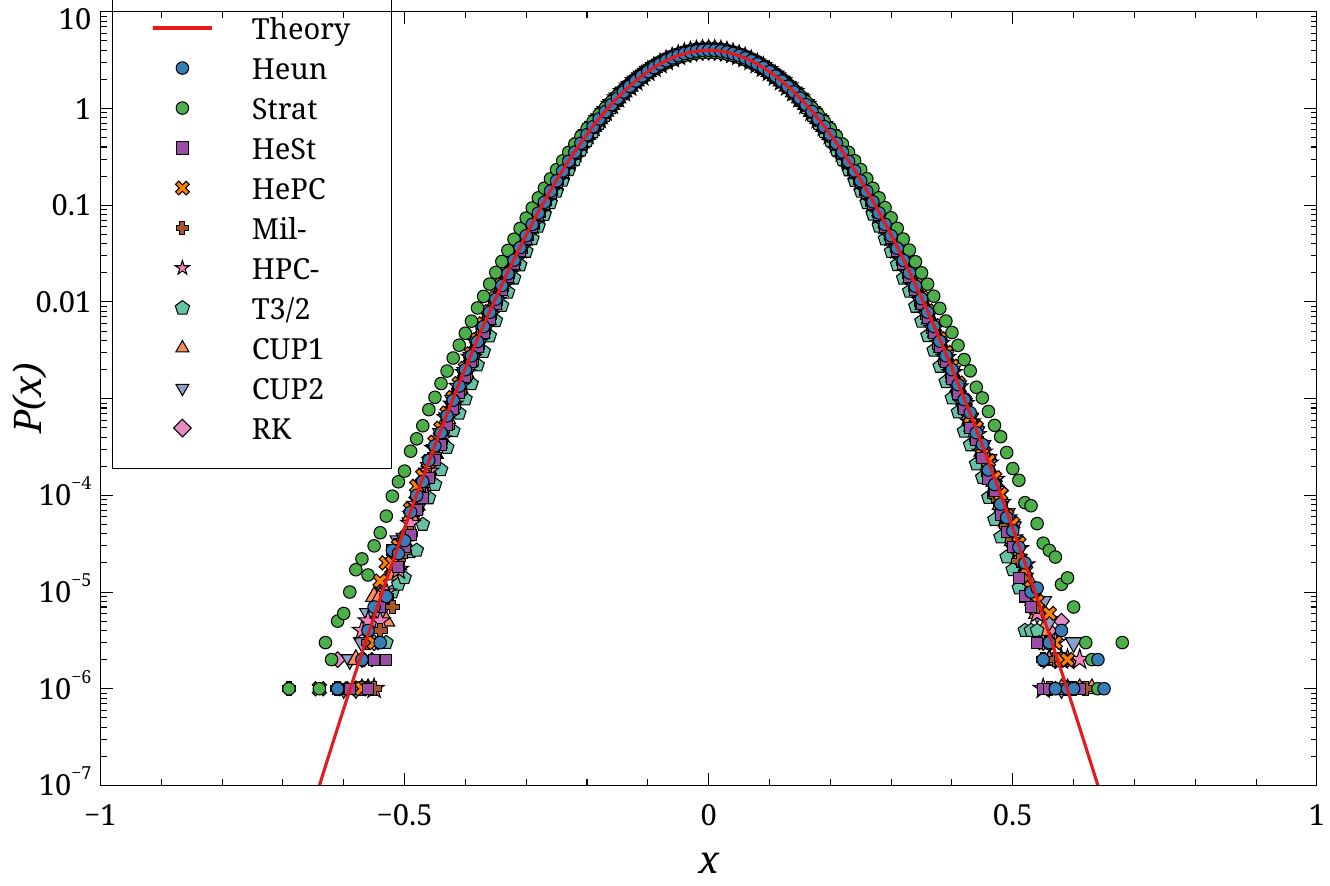} \\
\includegraphics[width=0.48 \textwidth]{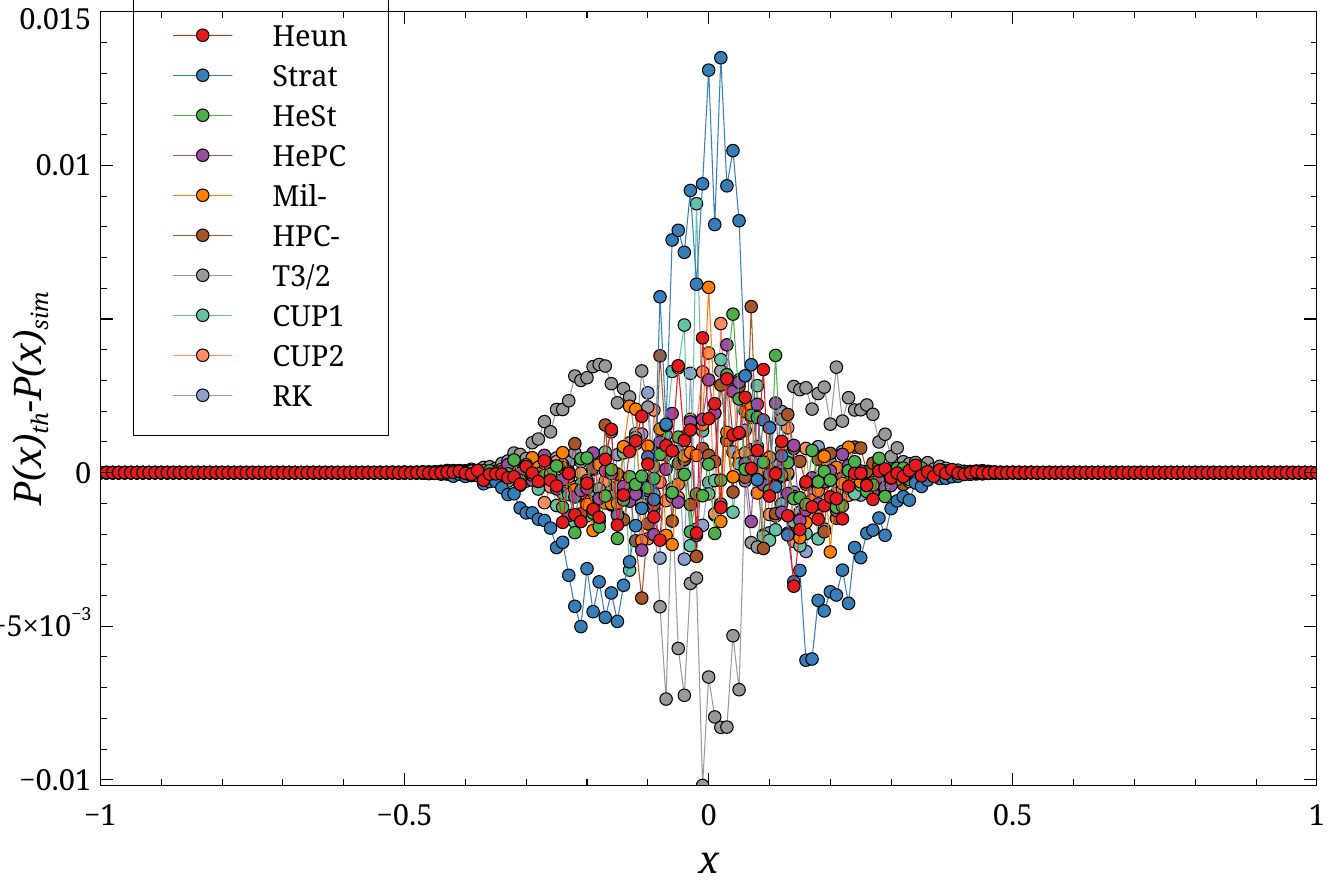}
\includegraphics[width=0.48 \textwidth]{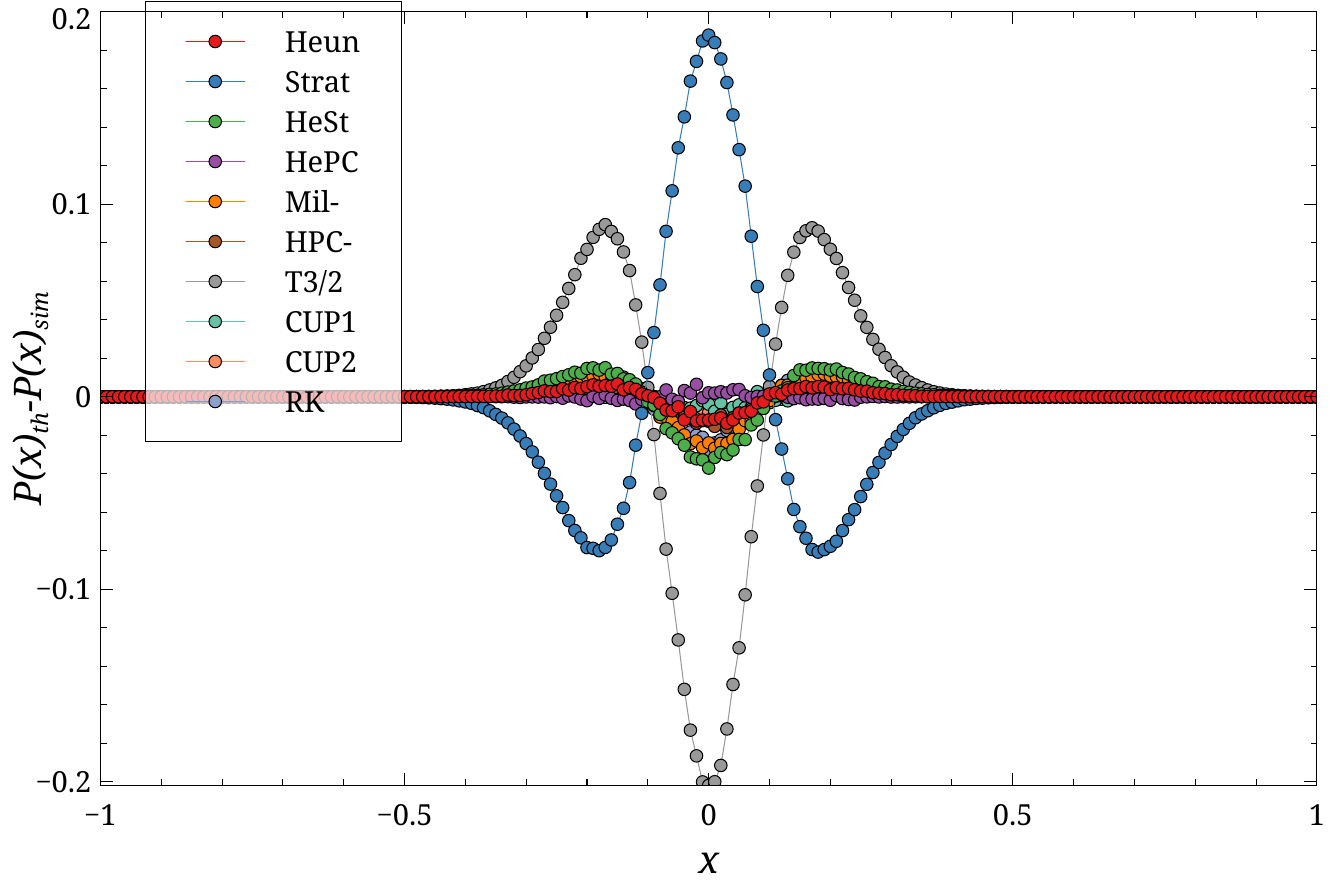}
\caption{Simulation results, done for $D=0.1 \times 10^{-1}$ and $h=0.01$ (left) and $h=0.2$ (right).
Top row shows $P(x)$, bottom row $P_{eq}(x)- P(x)$. Different symbols are the different schemes. \label{distributionsd}}
\end{figure}    

\begin{figure}[H]
\isPreprints{\centering}{} 
\includegraphics[width=0.48 \textwidth]{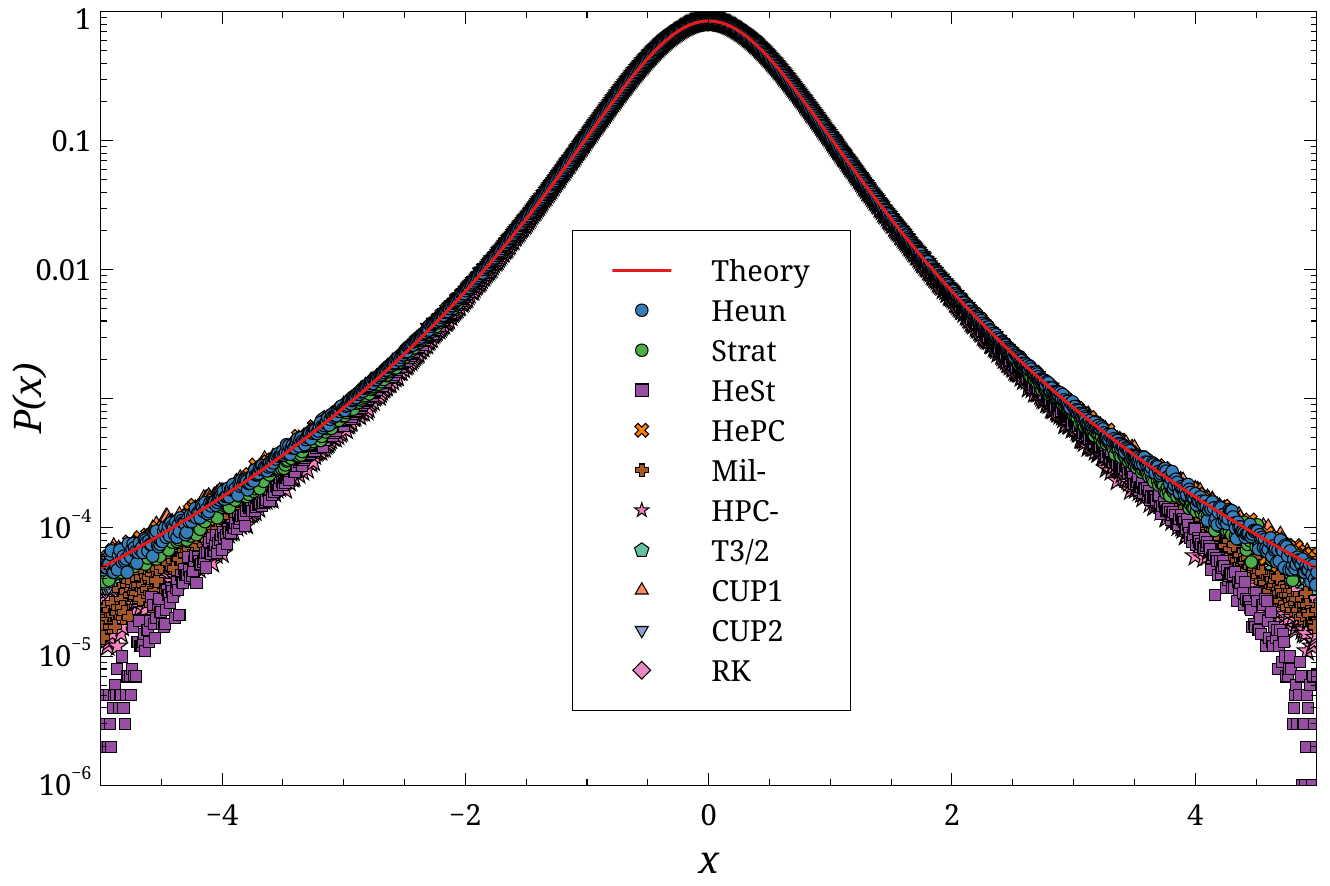}
\includegraphics[width=0.48 \textwidth]{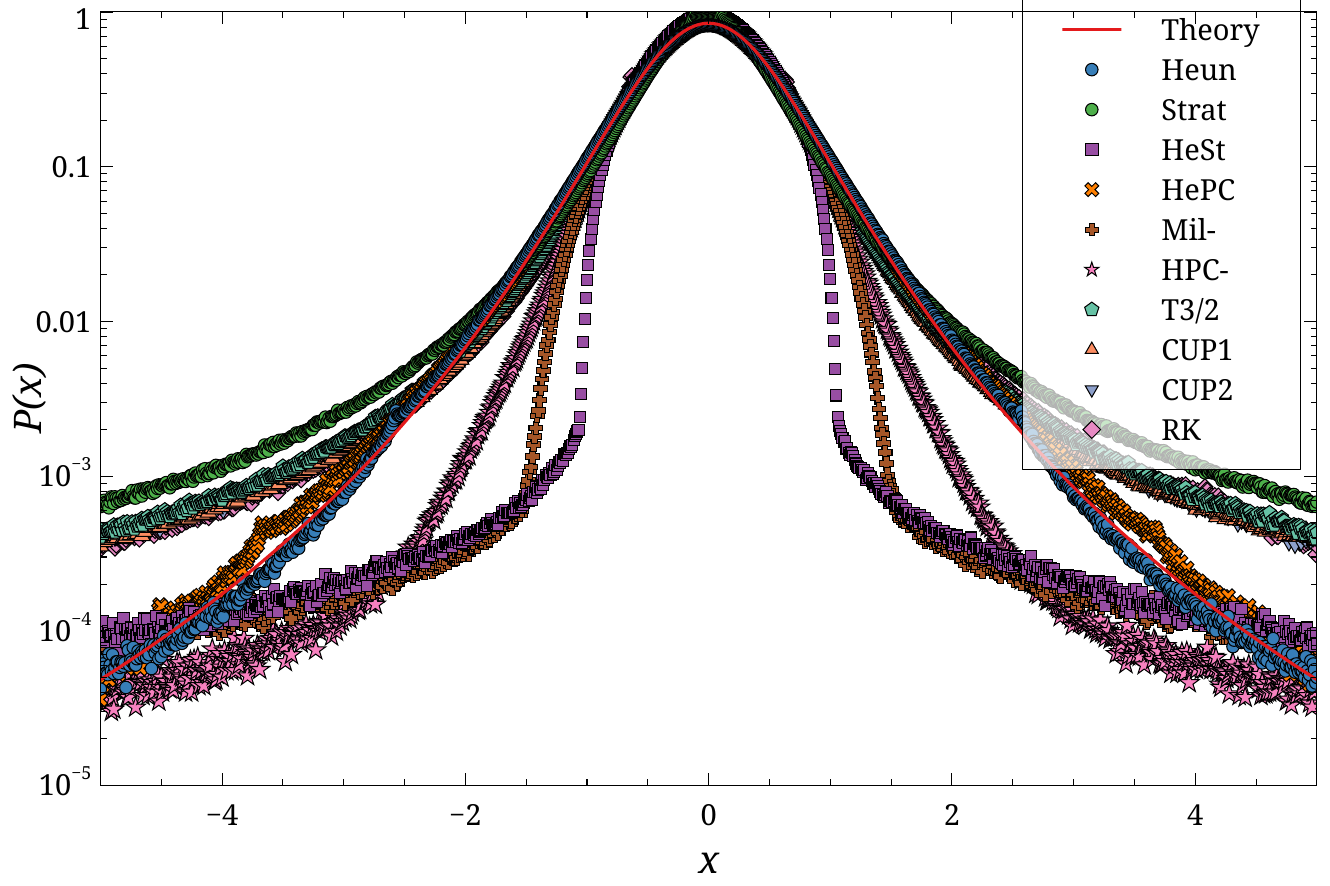} \\
\includegraphics[width=0.48 \textwidth]{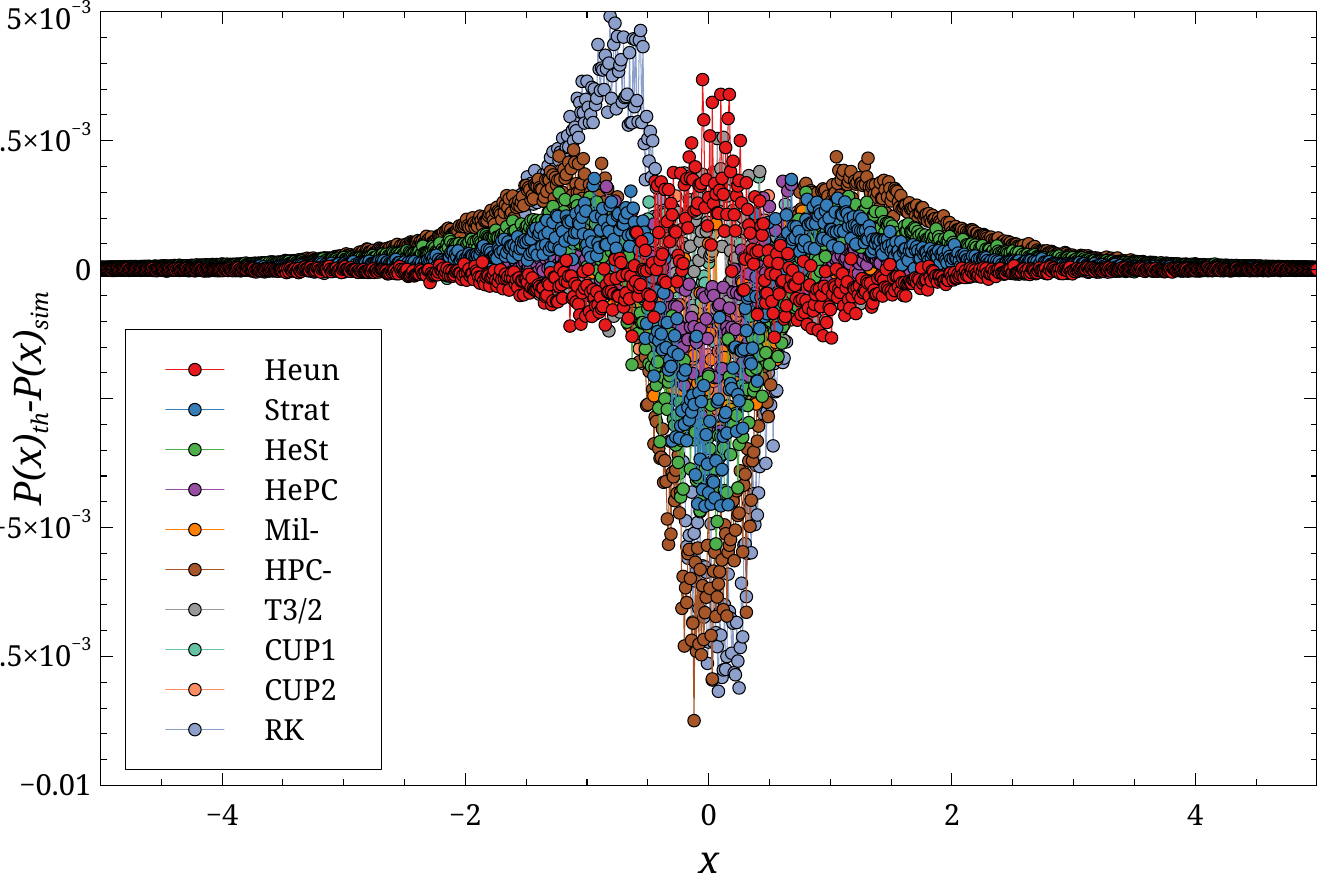}
\includegraphics[width=0.48 \textwidth]{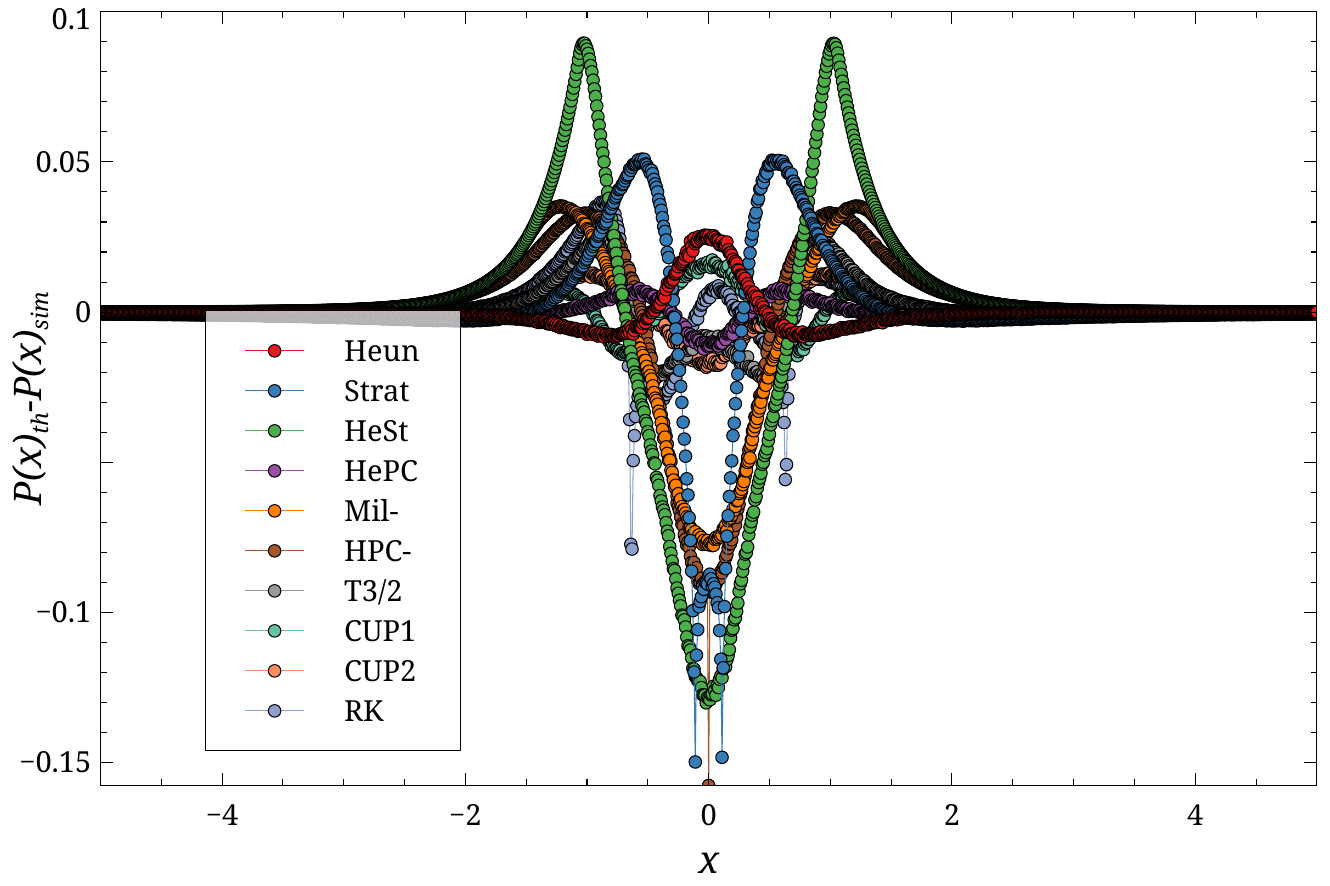}
\caption{Simulation results, done for $D=2.5 \times 10^{-1}$ and $h=0.01$ (left) and $h=0.2$ (right).
Top row shows $P(x)$, bottom row $P_{eq}(x)- P(x)$. Different symbols are the different schemes. \label{distributionld}}
\end{figure}    

Typical results are shown in Fig.~\ref{distributionsd} for a small $D$ value and in Fig.~\ref{distributionld} for a large $D$ value, and for two differente values of $h$.

These figures show that for small $h$, the agreement between theory and simulation is good for both values of $D$. As $h$ increases, the agreement deteriorates, with the magnitude of the deviation depending on the integration scheme. For instance, the results for $D=0.25$ and $h=0.2$ (Fig.~\ref{distributionld}, right)  show the emergence of spurious structures for some schemes. Overall, the \textbf{Stra} scheme is the first to deviate from the theoretical distribution. At $h=0.2$, the \textbf{HePC} scheme appears to yield the distribution closest to the theory, followed by \textbf{CUP2}, \textbf{CUP1}, and \textbf{Heun}, for both values of $D$.

To quantitatively assess the deviation of the simulated distribution $P(x)$ from the theoretical one $P_{eq}(x)$, we introduced two metrics. The first, which we term the "distance," is defined as $\int | P_{eq}(x)-P(x)| dx$. The second, the "ratio," is defined as $\int | (P_{eq}(x)-P(x))/P_{eq}(x)| dx = \int | 1-P(x)/P_{eq}(x)| dx$, where this integral is restricted to the region where $P_{eq}(x)\ge 10^{-4} P_{eq}(0)$. The distance metric emphasizes discrepancies near the origin where the probability density is large and statistics are robust, but it is less sensitive to differences in the tails. Conversely, the ratio metric places more weight on the tails, a region where statistics may be less reliable. For both metrics, a smaller value indicates better agreement between the simulated and theoretical distributions.

\begin{figure}[H]
\includegraphics[width=0.48 \textwidth]{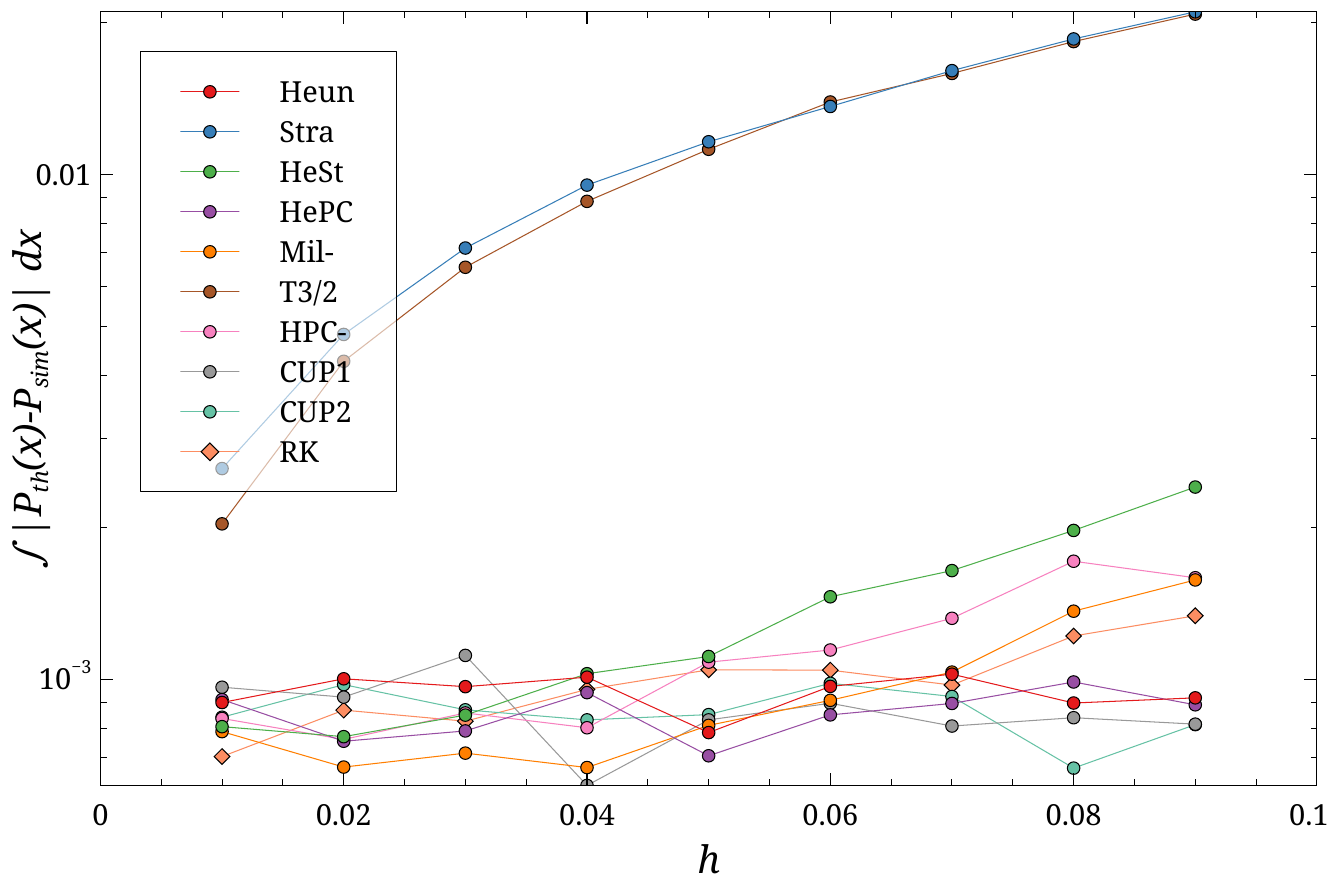}
\includegraphics[width=0.48 \textwidth]{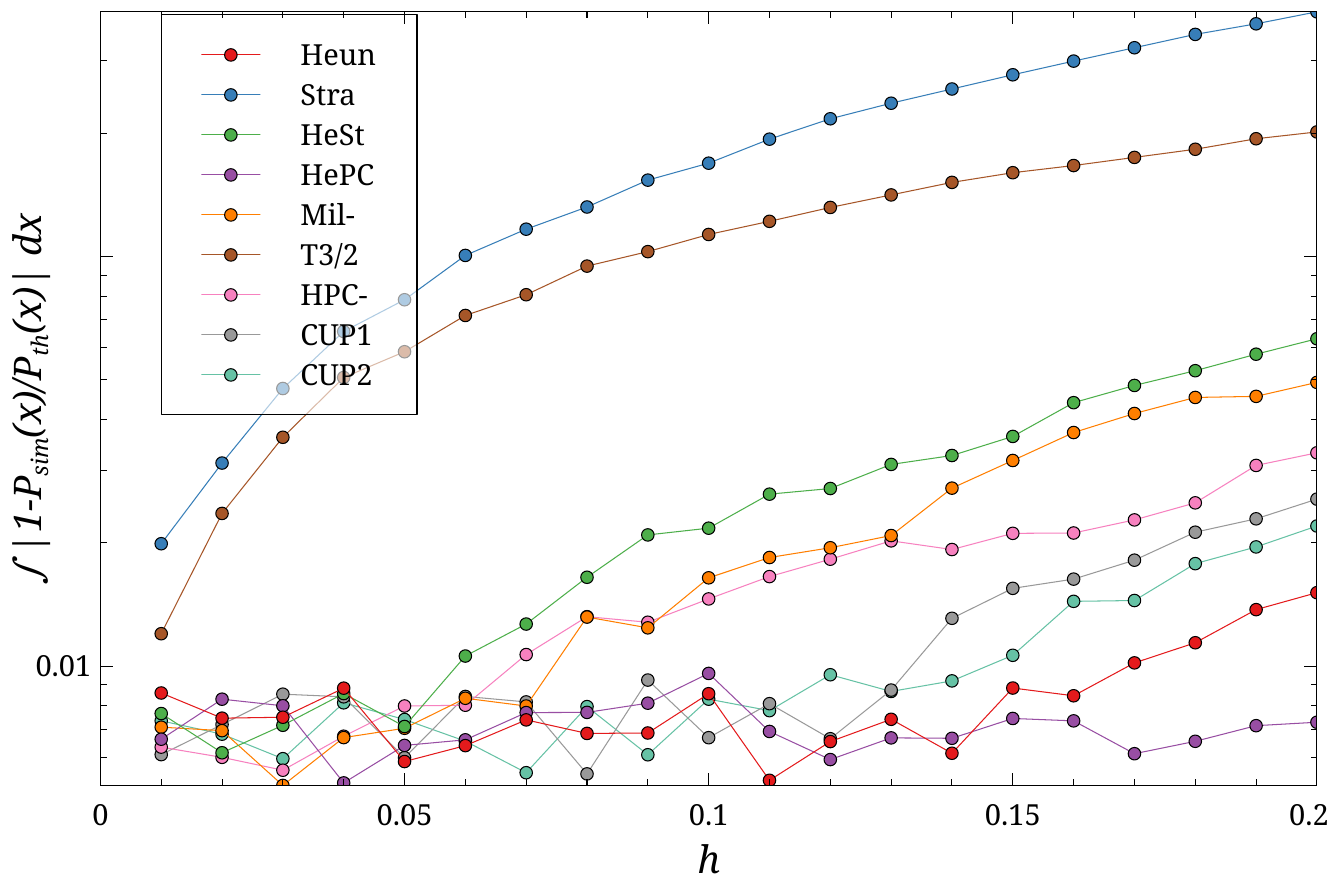}
\caption{The distance $\int | P_{eq}(x)-P(x)| dx$ (left) and the ratio $\int | 1-P(x)/P_{eq}(x)| dx$ (right) for
different $h$'s and $D=0.01$. Note the log scale on the ordinates. Different symbols are the different schemes. \label{disd0.01}}
\end{figure}

\begin{figure}[H]
\includegraphics[width=0.48 \textwidth]{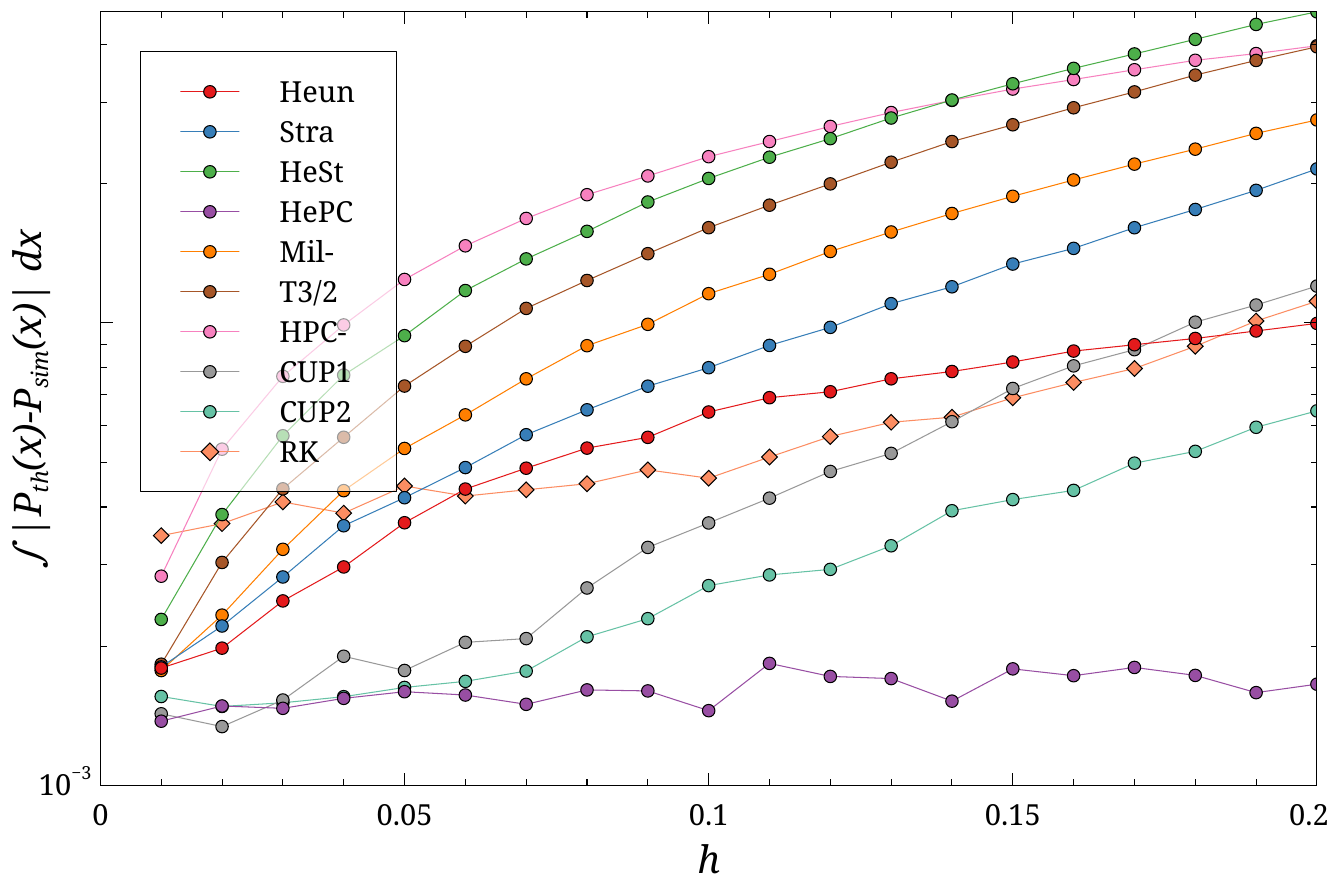}
\includegraphics[width=0.48 \textwidth]{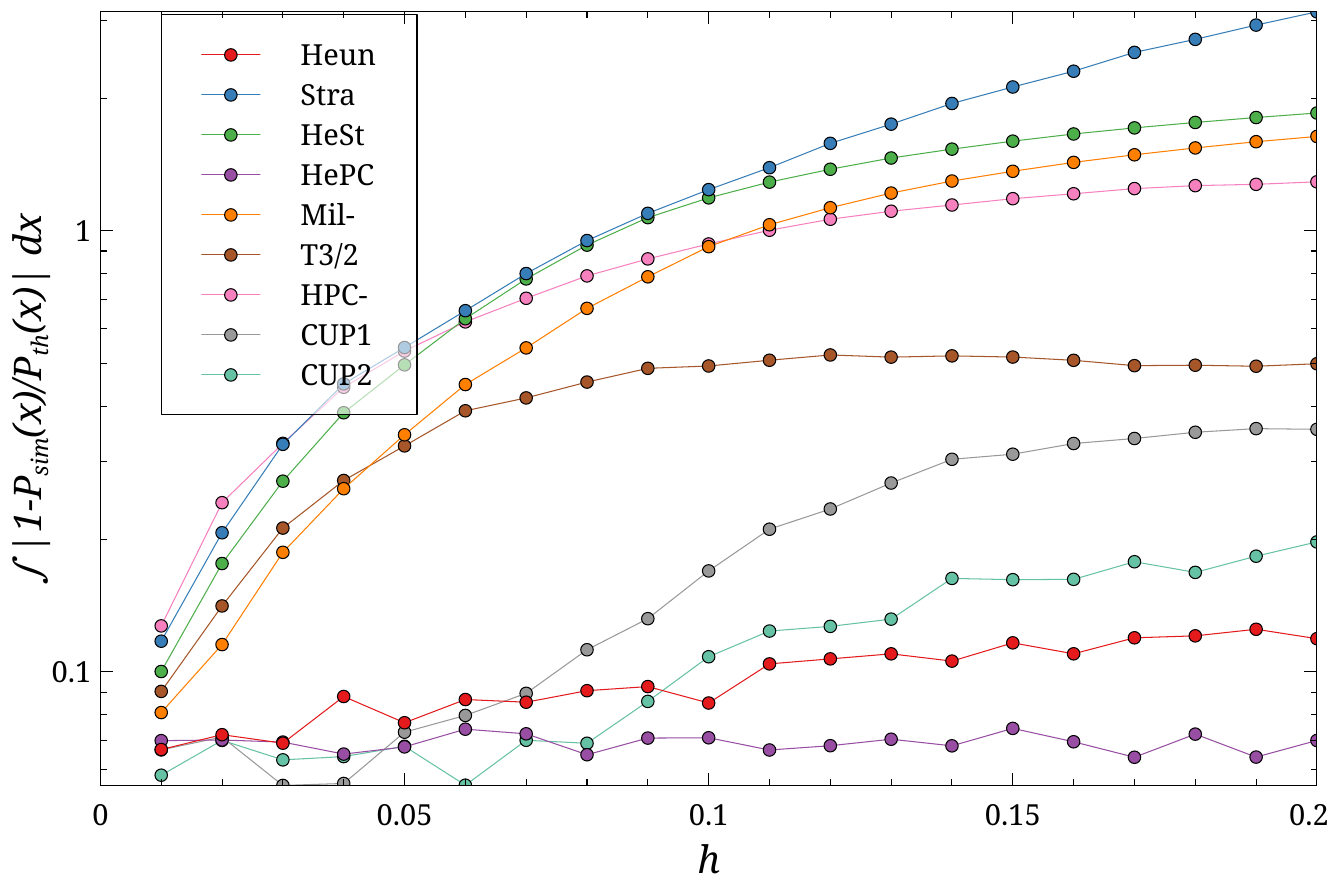}
\caption{The distance $\int | P_{eq}(x)-P(x)| dx$ (left) and the ratio $\int | 1-P(x)/P_{eq}(x)| dx$ (right) for
different $h$'s and $D=0.09$. Note the log scale on the ordinates. Different symbols are the different schemes. \label{disd0.09}}
\end{figure}

\begin{figure}[H]
\includegraphics[width=0.48 \textwidth]{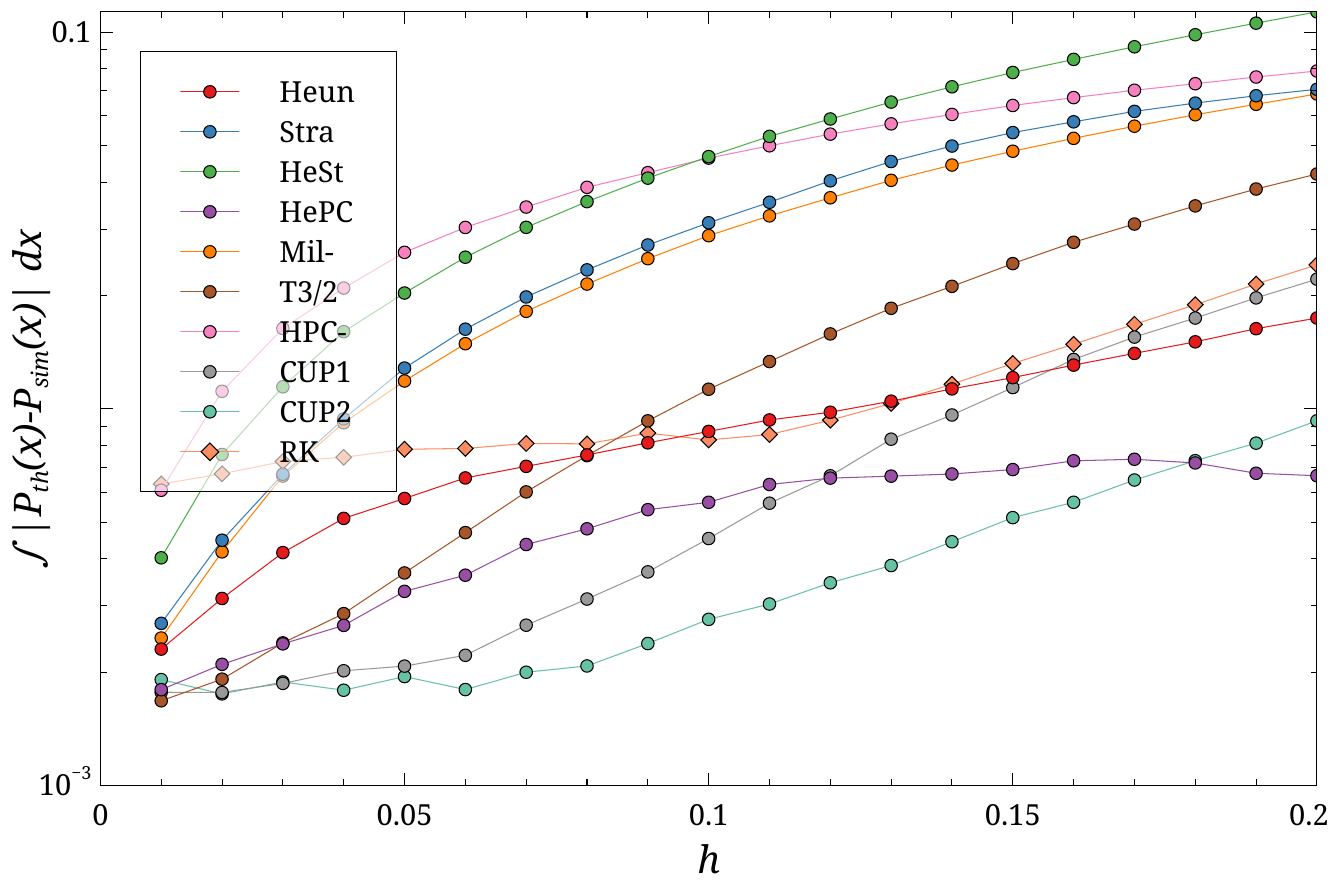}
\includegraphics[width=0.48 \textwidth]{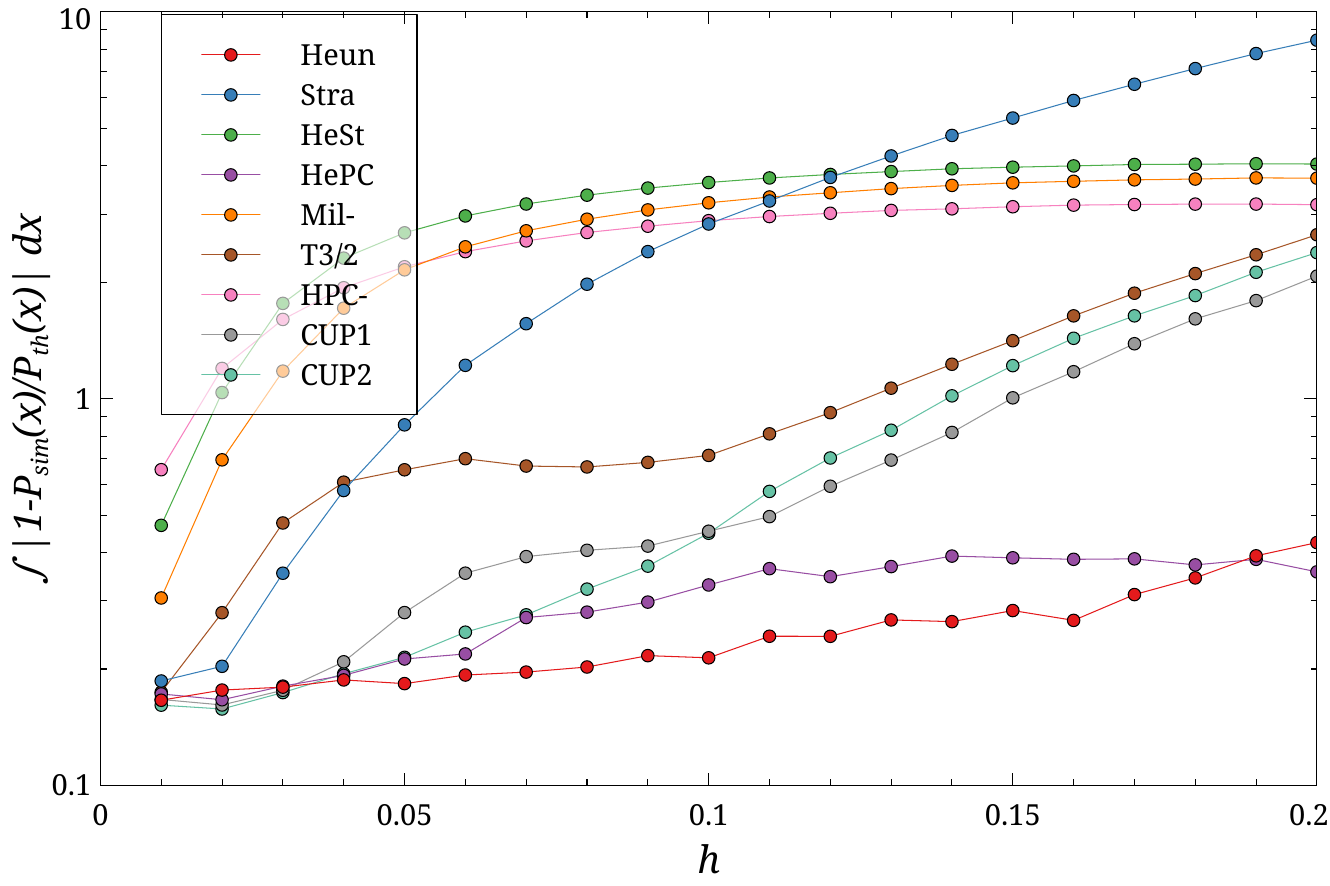}
\caption{The distance $\int | P_{eq}(x)-P(x)| dx$ (left) and the ratio $\int | 1-P(x)/P_{eq}(x)| dx$ (right) for
different $h$'s and $D=0.09$. Note the log scale on the ordinates. Different symbols are the different schemes. \label{disd0.17}}
\end{figure}

\begin{figure}[H]
\includegraphics[width=0.48 \textwidth]{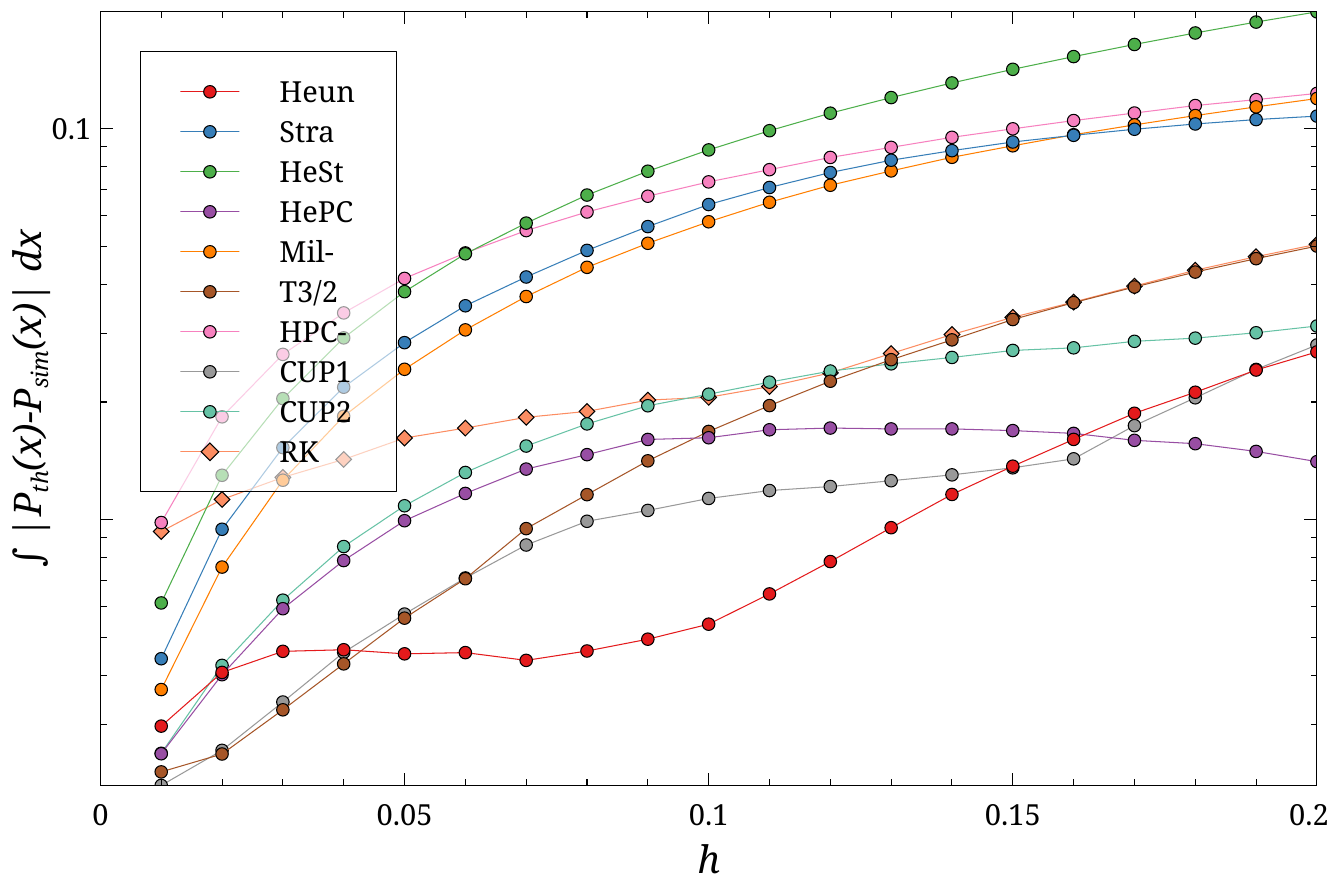}
\includegraphics[width=0.48 \textwidth]{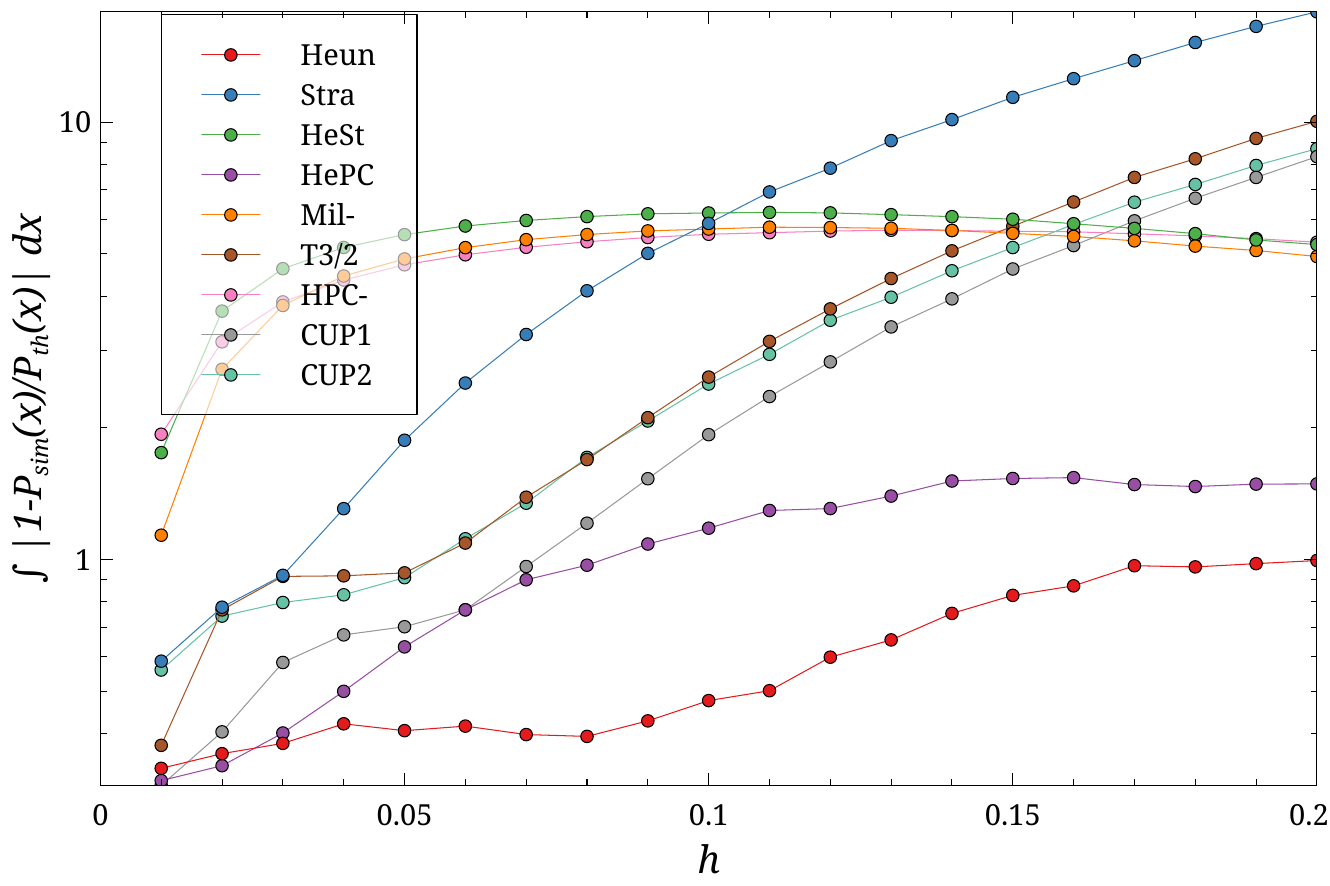}
\caption{The distance $\int | P_{eq}(x)-P(x)| dx$ (left) and the ratio $\int | 1-P(x)/P_{eq}(x)| dx$ (right) for
different $h$'s and$D=0.25$. Note the log scale on the ordinates. Different symbols are the different schemes. \label{disd0.25}}
\end{figure}

Figs~\ref{disd0.01},~\ref{disd0.09},~\ref{disd0.17} and~\ref{disd0.25} show the behaviour of the
two metrics as $h$ is increased, for four $D$ values.

As expected, both the distance and ratio metrics typically increase with $h$.
For the smallest value of $D$ considered (Fig.~\ref{disd0.01}), the distance between simulation and theory is small for most schemes at small $h$, with the \textbf{Stra} and \textbf{T3/2} schemes being notable exceptions. As $h$ increases, the distance for the \textbf{HePC} scheme remains roughly constant, while it grows for the others. At the largest $h$ considered, \textbf{HePC} performs best, followed by \textbf{CUP1}, \textbf{CUP2}, and \textbf{Heun}. The behavior of the ratio metric is similar, with only minor changes in the relative performance of \textbf{CUP1}, \textbf{CUP2}, and \textbf{Heun} at large $h$.

The evolution of the metrics for the largest $D$ considered (Fig.~\ref{disd0.25}) shows some interesting differences. At small $h$, the \textbf{HePC} scheme no longer stands out. It emerges as the best performer in terms of distance only as $h$ increases, and its advantage is less pronounced than for smaller $D$. In terms of the ratio, the \textbf{Heun} scheme performs best, followed by \textbf{HePC}.

The results for intermediate values of $D$, shown in Figs.~\ref{disd0.09} and~\ref{disd0.17}, exhibit a transitional behavior between these two extremes for both metrics.

\section{Conclusions}

In this work, we have conducted a critical re-examination of the Heun algorithm and several of its variants for the numerical integration of Stochastic Differential Equations. Motivated by the fact that the standard Heun scheme is based on the low-order Euler method, we systematically evaluated a dozen different algorithms, assessing their performance based on three key criteria: strong convergence, numerical stability, and the ability to reproduce correct long-time equilibrium distributions.

Our extensive numerical simulations lead to several key conclusions. First, despite its simplicity, the standard \textbf{Heun} scheme demonstrates remarkable robustness. It was found to be the most stable of all tested algorithms across a wide range of parameters and consistently produced accurate equilibrium distributions, especially when considering metrics that weight the tails of the distribution. This confirms its status as a reliable benchmark for SDE integration. There are some rare cases when \textbf{HePC} fared better, but given that this scheme requires several evaluation of the quantities $f(x)$ and $g(x)$ it should be weighted carefully whether to prefer \textbf{HePC}
to an integration done with \textbf{Heun} and a smaller $h$.

Second, a surprising and significant finding is the underperformance of theoretically higher-order schemes. Taylor-based schemes, such as \textbf{T3/2} and \textbf{CUP}, which include terms up to $O(h^{3/2})$ and $O(h^2)$ respectively, failed to achieve a strong convergence order beyond $O(h)$. This suggests that simply including more terms from a Taylor expansion does not guarantee improved strong convergence in practice, and highlights the complex relationship between theoretical order and numerical performance. These higher-order schemes also proved to be significantly less stable than the predictor-corrector type methods.

Third, our analysis of strong convergence has clarified a point of confusion in the literature. We have shown conclusively that the Heun scheme has a strong convergence of order $O(h)$, disproving an earlier claim that suggested an order of $O(h^{3/2})$.

\vspace{6pt} 



\funding{It is aknowledged support from the 
``National Centre for HPC,  Big Data and Quantum Computing'', under the National Recovery and 
Resilience Plan (PNRR), Mission 4 Component 2 Investment 1.4 funded from the 
European Union - NextGenerationEU.}

\dataavailability{Data are available upon reasonable requests. F90 codes are available at \url{https://github.com/dundacil/heun_extension}.} 

\acknowledgments{We thank the Green Data Center of University of Pisa for providing the computational  power needed for the present paper. During the preparation of this manuscript, the author used Gemini \url{https://gemini.google.com/app}  for the purposes of improving English and readability of some of the textual content. The author has reviewed and edited the output and take full responsibility for the content of this publication.}

\conflictsofinterest{The author declares no conflicts of interest.} 


\begin{adjustwidth}{-\extralength}{0cm}

\reftitle{References}


 \bibliography{heunpaper.bib}

\begin{thebibliography}{999}

\bibitem[Gardiner(2010)]{Gardiner_2010}
Gardiner, C.W.
\newblock {\em {Stochastic methods: A Handbook for the Natural and Social
  Sciences}}; Springer Series in Synergetics, Springer: Berlin, Germany,  2010.

\bibitem[Øksendal(2003)]{Oksendal_2003}
Øksendal, B.
\newblock {\em Stochastic Differential Equations. An introduction with
  applications}; Springer Berlin Heidelberg,  2003.
\newblock {\url{https://doi.org/10.1007/978-3-642-14394-6_5}}.

\bibitem[Kloeden and Platen(1992)]{Platen_Kloeden_1992}
Kloeden, P.E.; Platen, E.
\newblock {\em Numerical Solution of Stochastic Differential Equations};
  Springer Berlin Heidelberg,  1992.
\newblock {\url{https://doi.org/10.1007/978-3-662-12616-5}}.

\bibitem[Van~Kampen(2007)]{vanKampen_2007}
Van~Kampen, N.G.
\newblock {\em Stochastic Processes in Physics and Chemistry}; North Holland,
  2007.
\newblock {\url{https://doi.org/10.1016/B978-0-444-52965-7.X5000-4}}.

\bibitem[Mannella and McClintock(2012)]{Mannella_McClintock_2012}
Mannella, R.; McClintock, P.V.E.
\newblock It\^o versus Stratonovich: 30 years later.
\newblock {\em Fluctuation and Noise Letters} {\bf 2012}, {\em 11},~1240010,
  \href{http://arxiv.org/abs/https://doi.org/10.1142/S021947751240010X}{{\normalfont
  [https://doi.org/10.1142/S021947751240010X]}}.
\newblock {\url{https://doi.org/10.1142/S021947751240010X}}.

\bibitem[Heun(1900)]{HEUN1900}
Heun, K.
\newblock Neue Methode zur approximativen Integration der
  Differentialgleichungen einer unabh{\"a}ngigen Ver{\"a}nderlichen.
\newblock {\em Zeitschrift f{\"u}r Mathematik und Physik} {\bf 1900}, {\em
  45},~23--38.

\bibitem[Butcher(1996)]{BUTCHER1996247}
Butcher, J.
\newblock A history of Runge-Kutta methods.
\newblock {\em Applied Numerical Mathematics} {\bf 1996}, {\em 20},~247--260.
\newblock {\url{https://doi.org/https://doi.org/10.1016/0168-9274(95)00108-5}}.

\bibitem[R{\"u}melin(1982)]{Rumelin_1982}
R{\"u}melin, W.
\newblock Numerical treatment of stochastic differential equations.
\newblock {\em SIAM Journal on Numerical Analysis} {\bf 1982}, {\em
  19},~604–613.
\newblock {\url{https://doi.org/10.1137/0719041}}.

\bibitem[Kloeden and Pearson(1977)]{Kloeden_Pearson_1977}
Kloeden, P.E.; Pearson, R.A.
\newblock The numerical solution of stochastic differential equations.
\newblock {\em The Journal of the Australian Mathematical Society. Series B.
  Applied Mathematics} {\bf 1977}, {\em 20},~8–12.
\newblock {\url{https://doi.org/10.1017/S0334270000001405}}.

\bibitem[Klauder and Petersen(1985)]{Klauder_Petersen_1985}
Klauder, J.R.; Petersen, W.P.
\newblock Numerical Integration of Multiplicative-Noise Stochastic Differential
  Equations.
\newblock {\em SIAM Journal on Numerical Analysis} {\bf 1985}, {\em
  22},~1153--1166.

\bibitem[Mannella(2002)]{Mannella_2002}
Mannella, R.
\newblock Integration Of Stochastic Differential Equations On A Computer.
\newblock {\em International Journal of Modern Physics C (IJMPC)} {\bf 2002},
  {\em 13},~1177--1194.
\newblock {\url{https://doi.org/10.1142/S0129183102004042}}.

\bibitem[Mannella(1989)]{Mannella_1989}
Mannella, R., Computer experiments in non-linear stochastic physics.
\newblock In {\em Noise in Nonlinear Dynamical Systems}; Moss, F.; McClintock,
  P.V.E., Eds.; Cambridge University Press,  1989; p. 189–221.

\bibitem[Mil’shtejn(1975)]{Milshtejn_1975}
Mil’shtejn, G.N.
\newblock Approximate integration of stochastic differential equations.
\newblock {\em Theory of Probability \& Its Applications} {\bf 1975}, {\em
  19},~557–562.
\newblock {\url{https://doi.org/10.1137/1119062}}.

\bibitem[Rao et~al.(1974)Rao, Borwanker, and Ramkrishna]{Rao_1974}
Rao, N.J.; Borwanker, J.D.; Ramkrishna, D.
\newblock Numerical Solution of Ito Integral Equations.
\newblock {\em SIAM Journal on Control} {\bf 1974}, {\em 12},~124--139,
  \href{http://arxiv.org/abs/https://doi.org/10.1137/0312011}{{\normalfont
  [https://doi.org/10.1137/0312011]}}.
\newblock {\url{https://doi.org/10.1137/0312011}}.

\bibitem[Mannella and Palleschi(1989)]{Mannella1989FastAP}
Mannella, R.; Palleschi, V.
\newblock Fast and precise algorithm for computer simulation of stochastic
  differential equations.
\newblock {\em Physical Review A, General physics} {\bf 1989}, {\em
  40},~3381--3386.

\bibitem[Sancho et~al.(1982)Sancho, Miguel, Katz, and
  Gunton]{Sancho_Miguel_Katz_Gunton1982}
Sancho, J.M.; Miguel, M.S.; Katz, S.L.; Gunton, J.D.
\newblock Analytical and numerical studies of multiplicative noise.
\newblock {\em Physical Review A} {\bf 1982}, {\em 26},~1589--1609.
\newblock {\url{https://doi.org/10.1103/PhysRevA.26.1589}}.

\bibitem[Mannella(2004)]{Mannella_2004}
Mannella, R.
\newblock Quasisymplectic integrators for stochastic differential equations.
\newblock {\em Physical Review E} {\bf 2004}, {\em 69},~041107.
\newblock {\url{https://doi.org/10.1103/PhysRevE.69.041107}}.

\bibitem[Bogoi et~al.(2023)Bogoi, Dan, Strătilă, Cican, and
  Crunteanu]{Bogoi_et_al_2023}
Bogoi, A.; Dan, C.I.; Strătilă, S.; Cican, G.; Crunteanu, D.E.
\newblock Assessment of Stochastic Numerical Schemes for Stochastic
  Differential Equations with “White Noise” Using Itô’s Integral.
\newblock {\em Symmetry} {\bf 2023}, {\em 15}.
\newblock {\url{https://doi.org/10.3390/sym15112038}}.

\bibitem[mpi()]{mpich}
MPICH.
\newblock Accessed on Aug 10, 2025.

\bibitem[Press et~al.(1992)Press, Flannery, Teukolsky, and
  Vetterling]{press_etal:1992}
Press, W.H.; Flannery, B.P.; Teukolsky, S.A.; Vetterling, W.T.
\newblock {\em Numerical Recipes in FORTRAN 77: The Art of Scientific
  Computing}, 2 ed.; Cambridge University Press,  1992.

\bibitem[Higham(2001)]{Higham_2001}
Higham, D.J.
\newblock An Algorithmic Introduction to Numerical Simulation of Stochastic
  Differential Equations.
\newblock {\em SIAM Review} {\bf 2001}, {\em 43},~525--546,
  \href{http://arxiv.org/abs/https://doi.org/10.1137/S0036144500378302}{{\normalfont
  [https://doi.org/10.1137/S0036144500378302]}}.
\newblock {\url{https://doi.org/10.1137/S0036144500378302}}.

\bibitem[James and Roos(1975)]{James:1975dr}
James, F.; Roos, M.
\newblock {Minuit: A System for Function Minimization and Analysis of the
  Parameter Errors and Correlations}.
\newblock {\em Computer Physics Communications} {\bf 1975}, {\em 10},~343--367.
\newblock {\url{https://doi.org/10.1016/0010-4655(75)90039-9}}.

\bibitem[Honeycutt(1992)]{Honeycutt_1992}
Honeycutt, R.L.
\newblock Stochastic Runge-Kutta algorithms. I. White noise.
\newblock {\em Physical Review A} {\bf 1992}, {\em 45},~600--603.
\newblock {\url{https://doi.org/10.1103/PhysRevA.45.600}}.

\end{thebibliography}

\PublishersNote{}
\end{adjustwidth}
\end{document}